\documentclass[11pt,twoside,a4paper]{article}
\usepackage{garamondx}
\usepackage[polutonikogreek ,british]{babel}
\usepackage{soul}
\usepackage{multicol}
\usepackage{amsmath,amssymb}
\usepackage{comment}
\usepackage{cancel}
\usepackage{flushend}
\usepackage{fancybox}
\usepackage[utf8]{inputenc}
\usepackage{graphics,epsfig} 
\usepackage{floatflt}
\usepackage{float}
\usepackage{longtable}
\usepackage{fancybox}
\usepackage{subfig}
\usepackage[pdftex]{color}
\usepackage{yhmath}
\usepackage{xcolor}
\usepackage{tikz}
\usepackage{tkz-euclide}

\definecolor{urlcolor}{rgb}{0,0,0}
\definecolor{linkcolor}{rgb}{0,0,0}
\definecolor{citecolor}{rgb}{0,0,0}
\usepackage[pdfpagemode=UseOutlines,urlcolor=urlcolor,linkcolor=linkcolor,citecolor=citecolor,colorlinks=true,backref=page]{hyperref}
\usepackage[nameinlink]{cleveref}  
\renewcommand*{\backref}[1]{}

\usepackage{paracol}

\makeatletter
\def\section{\@startsection {section}{1}{\z@}{-3.5ex plus -1ex minus
-.2ex}{2.3ex plus .2ex}{\large\bf}}
\def\subsection{\@startsection {subsection}{2}{\z@}{-3.25ex plus -1ex minus
-.2ex}{1.5ex plus .2ex}{\bf}}

\addto\captionsgerman{} 

\begin{document}
 
\thispagestyle{empty}
\title{\Large      \bf
{
\ovalbox{\ovalbox{
\begin{minipage}{9cm}
\begin{center}
~\\[1.5ex]
Benedetti's constructions of ovals\\
~\\
\end{center}
\end{minipage}
}\!\!
}}}
\author{\normalsize  
\textbf{Thomas Hotz}\thanks{\tt thomas.hotz{\scriptsize @}tu-ilmenau.de}
 \quad\text{and}   \quad
 \textbf{Achim Ilchmann}\thanks{\tt achim.ilchmann{\scriptsize @}tu-ilmenau.de, corresponding author
}
 \\[0.5ex]
  Institut f\"ur Mathematik, Technische Universit\"{a}t Ilmenau \\[0.5ex]
   Weimarer Stra{\ss}e~25, 98693 Ilmenau \\[-1ex]
}
\date{(submitted to \textit{Historia Mathematica}, 12~December~2024)}

\maketitle
\vspace*{-1.3cm}
\tableofcontents
\thispagestyle{empty}

\newpage
\section*{Abstract}  
Giovanni Battista Benedetti~(1530--1590) derived
two constructions of ovals given their minor and major axes. These were published in~1585 and seem to be the first solution to this problem. Therefore, the generally accepted view that   ``the geometrical construction for drawing an oval for any given
proportion was not known in the sixteenth century''~(Ana López-Mozo, 2011; Nexus Network Journal, p.\,571)
has to be revised.
We provide a transcription of Benedetti's Latin text as well as a translation into English and analyze his mathematical reasoning in detail. 
This is the first main focus of the present note. The second focus is the impact of Benedetti's contributions to ovals. 
In particular, it is observed that they were rarely mentioned in the mathematical and architectural literature.

\setcounter{page}{2}

\section{A brief introduction to ovals}\label{Sec:intro}  
``The oval shape is a feast for the eyes, perfect, and due to its simplicity and softness, it can be used for various purposes.''
concludes Sebastiano Serlio (1475~-- ca.\,1554) in his 
\textit{Il primo libro d'architettura} from~1545.~\cite[fol.\,~20v]{Serlio-1545-engl}
Completeness and simplicity  are a necessity for beauty,
which is realised through geometric constructions such as oval forms.
These forms played an important r\^{o}le in Baroque architecture.~\cite{Chatelet-Lange-1976-Engl,Lotz-1955-engl}
They were publicised in
Serlio's  treatise \textit{Sette Libri d'Architettura} 
which was, in Santiago Huerta's words, 
   ``one of the most popular architectural treatises ever published''.~\cite[p.\,230]{Huerta-2007}
Numerous editions and receptions of it were published in Baroque times.~\cite[p.\,215]{Oechslin-Buechi-Pozsgai-2018}, \cite[pp.\,80--87]{Kruft-1985-engl}

\begin{figure*}[t!]
\begin{center}
\subfloat[Oval~S1: two equilateral  triangles]
 {
 \centering
   \includegraphics[angle=0,width=6.3cm]{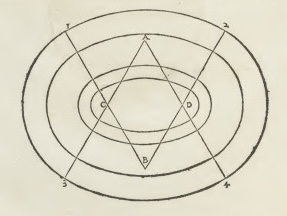}
\label{Abb:Serlio1}
 \vspace{5mm}
 }
 \hspace*{12mm}
\subfloat[Oval~S2: two adjacent equal circles]
 {
  \centering
   \includegraphics[angle=0,width=6.3cm]{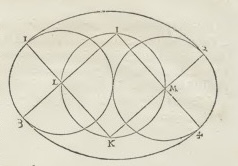}
 \label{Fig:Serlio2}
}
 \qquad\quad
\subfloat[Oval~S3: two adjacent equal squares]
 {
  \centering
   \includegraphics[angle=0,width=6.3cm]{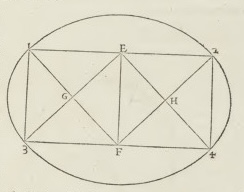}
\label{Abb:Serlio3}
 }
  \hspace*{12mm}
\subfloat[Oval~S4: two intersecting equal circles]
 {
  \centering
   \includegraphics[angle=0,width=6.3cm]{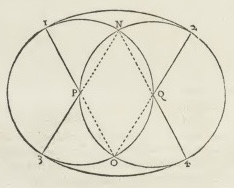}
\label{Abb:Serlio4}
 }
 \caption{\small  Serlio's four 4-centre ovals with right-angled, isosceles~(S2, S3) and
 equilateral~(S1, S4) triangles; details from~\cite[pp.\,17--19]{Serlio-1545-engl}.}
 \label{Fig:Serlio-ovals}
 \end{center}
\end{figure*}

However,  at those times~-- and even today~--  there has not been a
consistent  definition of an oval;  
a mathematician, an architect and an art historian would each  
suggest  a different definition. In the present note we consider only what is nowadays called a \textit{4-centre oval},
which is a smooth, convex, closed curve formed by four circular arcs being symmetric with respect to two perpendicular axes. The \textit{centres} of an oval, which  form a rhombus, are those points that serve as the centres of the circular arcs constituting the oval. 
Here, smoothness means that the circles touch each other at the connection points of the respective arcs, i.e., their tangent lines at that point agree; equivalently, 
as Euclid showed~\cite[Prop.\,III.11]{Euclid-Clavius-1591-engl},
the circles' centres and their connection point are collinear.~(Fig.~\ref{Fig:Serlio-ovals})

We investigate~-- and this is the core of the present note~-- 
the contributions to ovals by the Italian mathematician Giovanni Battista Benedetti (1530–1590), aa well as the treatment of his findings in architectural and mathematical publications. Benedetti discovered how to construct an oval 
when  two semiaxes and the radius of the smaller (or alternatively the larger) circular arcs are given; we will call the first one \textit{oval~B1} and the second one \textit{oval~B2}. The novelty of his findings becomes obvious when compared to Serlio's ovals, as no other ovals were known at that time. Ana López Mozo writes  that ``none of Serlio's ovals can be adjusted to fit any given place,
but only those which have the same proportions''~\cite[p.\,572]{Lopez-Mozo-2011}. Benedetti, in fact, solved this architecturally important problem.

Benedetti discovered his oval constructions while solving another problem in Euclidean geometry.~(Section~\ref{Sec:constructions})
There are no other known results by Benedetti on ovals. 
For him, it seems to have been an intriguing mathematical problem in its own right, he himself being
unaware of its significance to architecture or mathematics. 
He reported his findings in a letter to the nobleman Pietro Pizzamano (1512--1571), which he  must have written before Pizzamano's death in~1571. It was published later in 1585, however, buried in his last work~\cite{Benedetti-1585} under a heap of other miscellanea, so it is no wonder it went largely unnoticed. (Section~\ref{Sec:Reception}) 
  
As a consequence, Ana López Mozo's  conjecture,  
``The Spanish Tomás Vicente Tosca~[1712] might possibly have 
given the first, general construction for drawing an
oval for any given proportion.''~\cite[p.\,580]{Lopez-Mozo-2011}
 must be corrected by at least~141 years.

To present an (incomplete)  list of   references  on ovals, we mention that many historical ovals, as well as systematically ordered ovals, are discussed in Angelo Mazzotti's monograph \textit{All Sides to an Oval}.~\cite{Mazzotti-2019} Ovals from the sixteenth century onwards are studied by Ana López Mozo~\cite{Lopez-Mozo-2011}, and a source on oval Roman amphitheatres is Mark Wilson Jones' article~\cite{Wilson-Jones-1993}. This is not an exhaustive list of the literature on the subject, of course; further references are cited in the works mentioned above. 
 
\captionsetup[subfloat]{labelformat=empty}

\section{Giovanni Battista Benedetti and his letters to Pietro Pizzamano (before~1571)}
\label{sec:benedetti}
%
 

Giovanni Battista Benedetti, also known as Giambattista Benedetti, was born on 14~August~1530 in Venice and died on 20~January~1590 in Turin.\footnote{For bibliographical details on Benedetti's life, we refer to the established literature, in particular~\cite{Bauer-1991,Bordiga-1926,Omodeo-Renn-2019,Roero-1997}.} 
In his first work, which appeared in~1553 at Venice, he demonstrated that all of Euclid's straightedge-and-compass construction problems can also be solved when using a compass with a fixed opening.
This established his reputation as an excellent geometer.\footnote{In retrospect, Moritz Cantor called him a ``true geometer''.~\cite[p.\,521]{Cantor-1892}}
In~1558, he accepted a position as lecturer of mathematics and philosophy at the Farnese court in Parma, from where he moved to Turin in~1567 to become the court mathematician of the Duke of Savoy, a position he held until his death.

Throughout his life, Benedetti worked on various topics, ranging from arithmetics to music.
In fact, he is well-known today for his contributions not only to geometry, but also to perspective~\cite{Andersen-2007,Field-1985}  and theoretical mechanics~\cite{Bauer-1991,Omodeo-Renn-2019}.
This diversity of research interests is reflected in the title of his last work, published first in~1585 at Turin, called \textit{Diversarum speculationum mathematicarum et physicarum liber}~\cite{Benedetti-1585}, or \textit{Diversarum speculationum} for short; two further editions appeared in Venice in~1586 and~1599 under slightly varying titles~\cite[p.~28]{Omodeo-Renn-2019}.
It comprises five treatises (on arithmetics, perspective, mechanics, proportions, and on some of Aristotle's views) along with~118 letters on a multitude of topics.\footnote{More information on the \textit{Diversarum speculationum} may be found in the literature~\cite{Bauer-1991,Bordiga-1926,Omodeo-Renn-2019}.}

Here, we will be concerned with just two of these letters.~\cite[pp.\,262--264]{Benedetti-1585}
Both are addressed to \textit{Pietro di Giovanni Andrea Pizzamano} (1512--1571)\footnote{See \cite[vol. 6=ms.~2503, fols.\ 103v--104r]{Barbaro} for Pizzamano's genealogy including dates of birth and death.} who was a patrician of the Republic of Venice.
He served as Savio agli ordini (a magistrate for maritine matters), as Savio alla mercancia (magistrate for commerce), as Provveditore del sal (superintendent for salt), as Console della Soria (councillor), as Capitano of Bergamo, as Podestà as well as Capitano of Bassano and of Trevigi, and as Rector of Sitia 
(on Crete).~\cite{Bordiga-1926,Cecchini-Roero-2004,Szepe-2018}
We note that, as Pizzamano died in~1571, it can safely be assumed the letters to him have been written no later than that, i.e.\ at least~14 years before the first edition of the \textit{Diversarum speculationum} appeared. 

The interested reader may find two miniatures portraying Pizzamano in \cite{Szepe-2018}, whereas the only portrait of Benedetti was lost to a fire at the national library at Turin in~1904 \cite[p.~603]{Bordiga-1926}, as were two large collections of his letters~\cite[p.~28]{Omodeo-Renn-2019}.






\section{Benedetti's oval and encompassing circle constructions}
\label{Sec:constructions}
The two letters in Benedetti's \textit{Diversarum Speculationum} mentioned above are concerned with the construction of certain circles and ovals, see~\cite[pp.\,262--264]{Benedetti-1585}
and their transcription as well as their translation into English in the Appendix.
In the first letter, he introduces a construction  for a circle
encompassing two given circles, 
which is treated in Section~\ref{Ssec:Bcircles}.
In the second letter,  he presents
two general constructions of ovals where both the major and minor semiaxes, as well as the radius of the smaller or larger circular 
arcs are given;
these are treated in Section~\ref{Ssec:B1} and~\ref{Ssec:B2}.
The constructions of ovals~B1 and~B2  in the second letter stand on their own and are more important than the results in the first letter. 
 However, as we will see, the problem addressed in the first letter is closely related to Benedetti's construction of ovals.

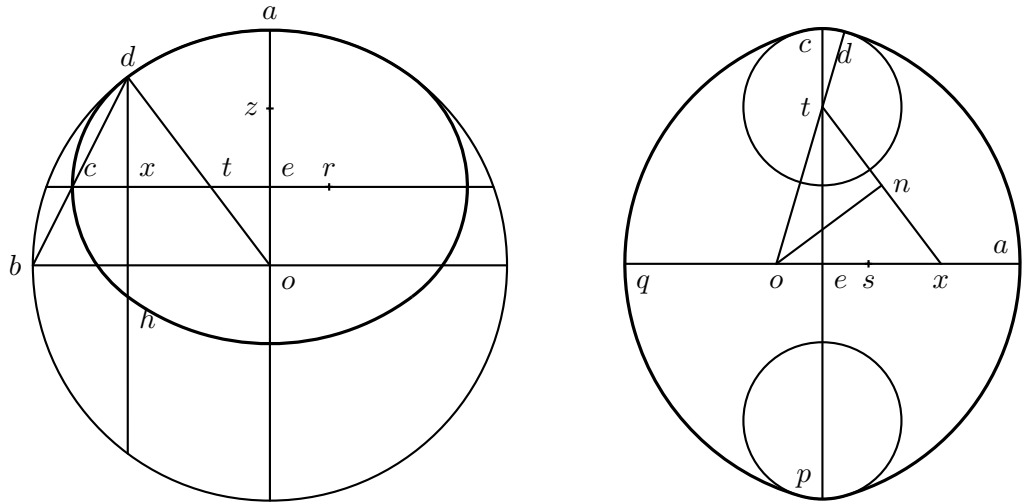
\begin{figure*}[t!]
   \begin{center}
 \subfloat[(a)  \ Oval B1: given semiaxes~$ea<ec$ and the radius~$oa$ of  the larger circular arc]{
 \centering
 \begin{tikzpicture}[scale=2.6]
  \tkzDefPoints{0/0/e, 0/0.8/a, -1/0/c, 0/-0.4/o, -1.2/-0.4/b, 1.2/-0.4/B, 1.2/-0.4/b', 1/0/c'}
  \tkzDrawCircle[black, thick](o,a)
  \tkzDrawSegment[thick](b,B)
  \tkzInterLC(o,a)(o,a) \tkzGetFirstPoint{A}
  \tkzDrawSegment[thick](a,A)
  \tkzInterLC(e,c)(o,a) \tkzGetPoints{cc}{C}
  \tkzDrawSegment[thick](cc,C)
  \tkzInterLC(b,c)(o,a) \tkzGetFirstPoint{d}
    \tkzInterLC(b',c')(o,a) \tkzGetSecondPoint{d'} 
  \tkzDrawSegments[thick](b,d d,o)
  \tkzInterLL(d,o)(c,e) \tkzGetPoint{t}
  \tkzDefLine[perpendicular=through d](e,c) \tkzGetPoint{hd}
  \tkzInterLL(d,hd)(c,e) \tkzGetPoint{x}
  \tkzInterLC(d,hd)(t,d) \tkzGetSecondPoint{h}
  \tkzDrawArc[color=black, very thick](t,d)(h)
    \tkzDrawArc[color=black, very thick](o,d')(d)   
  \tkzInterLC(d,hd)(o,a) \tkzGetSecondPoint{D}
  \tkzDrawSegment[thick](d,D)
  \tkzInterLC(a,e)(e,o) \tkzGetFirstPoint{z}
  \tkzDrawPoint[shape=cross, thick](z)
  \tkzInterLC(h,e)(o,a) \tkzGetFirstPoint{hh}
  \tkzInterLC(c,e)(e,t) \tkzGetSecondPoint{r}
  \tkzDrawPoint[shape=cross, thick](r)
  \tkzInterLC(d,e)(z,h) \tkzGetSecondPoint{dd}
  \tkzDrawArc[color=black, very thick](z,h)(dd)
  \tkzDrawArc[color=black, very thick](r,dd)(hh)
  \tkzLabelPoints[above](a,d,r)
  \tkzLabelPoints[left](b,z)
  \tkzLabelPoints[below right](o,h)
  \tkzLabelPoints[above right](e,c,t,x)
\end{tikzpicture}
\label{Abb:Benedetti-B1}
 \vspace{4mm}
}
 \qquad 
 \subfloat[(b) \ Oval B2: given semiaxes~$ea<ec$ and the radius~$ct$ of the smaller circular arc]{
 \centering
\begin{tikzpicture}[scale=2.6]
  \tkzDefPoints{0/0/e, 0/1.2/c, 1/0/a, 0/0.8/t, 0/-0.8/t'}
  \tkzInterLC(e,c)(e,c) \tkzGetFirstPoint{p}
  \tkzInterLC(e,a)(e,a) \tkzGetFirstPoint{q}
  \tkzDrawSegments[thick](c,p a,q)
  \tkzDrawCircle[black, thick](t,c)
  \tkzCalcLength(t,c) \tkzGetLength{tc}
  \tkzInterLC[R](e,a)(a,\tc) \tkzGetFirstPoint{x}
  \tkzDrawSegment[thick](t,x)
  \tkzDefMidPoint(t,x) \tkzGetPoint{n}
  \tkzDefLine[orthogonal=through n](t,x) \tkzGetPoint{no}
  \tkzInterLL(n,no)(a,e) \tkzGetPoint{o}
  \tkzDrawSegment[thick](n,o)
  \tkzInterLC[near](o,t)(t,c) \tkzGetSecondPoint{d} 
     \tkzDrawSegment[thick](o,d)
  \tkzInterLC(e,o)(e,o) \tkzGetFirstPoint{s}
  \tkzDrawPoint[cross, thick](s)
  \tkzInterLC(e,t)(e,t) \tkzGetFirstPoint{T}
  \tkzDrawCircle[black, thick](T,p)
  \tkzInterLC[near](o,T)(T,p) \tkzGetSecondPoint{D} 
\tkzDrawArc[color=black, very thick](o,D)(d)
  \tkzInterLC(e,d)(T,p) \tkzGetSecondPoint{dd}
  \tkzInterLC(e,D)(t,c) \tkzGetFirstPoint{DD}
  \tkzDrawArc[color=black, very thick](s,DD)(dd) 
  \tkzLabelPoints[below](x,o,d,s)
  \tkzLabelPoints[left](t)
  \tkzLabelPoints[right](n)
  \tkzLabelPoints[above left](a,p)
  \tkzLabelPoints[below left](c)
  \tkzLabelPoints[below right](e,q)

 \tkzInterLC[near](s,t)(t,c) \tkzGetSecondPoint{d'}  %
   \tkzDrawArc[color=black, very thick](t,d)(d')    %
   
     \tkzInterLC(o,t')(t',p) \tkzGetFirstPoint{f}
      \tkzInterLC(s,t')(t',p) \tkzGetSecondPoint{f'}
       \tkzDrawArc[color=black, very thick](t',f')(f)    %
\end{tikzpicture}
\label{Abb:Benedetti-B2}
 } 
 \vspace{1.1ex}
 \caption{\small  Benedetti's construction of ovals~B1 and~B2 (recreated; originals in Fig.~\ref{Fig:Benedetti-264-B1-B2}).}
\label{Abb:Benedetti-Ovals}
 \end{center}
\end{figure*}

\subsection[Benedetti's oval~B1]{Benedetti's oval~B1\quad \normalfont{(Fig.\,\ref{Abb:Benedetti-B1})}}\label{Ssec:B1}

 In the first part of the second letter, 
 Benedetti constructs  an oval  
given both the major and minor semiaxes as well as
the length of the  radii of the  larger circular arcs.~(Fig.\,\ref{Abb:Benedetti-B1})
We call this construction \textit{oval~B1}.
Specifying these three parameters determines the oval uniquely, as will be seen.
 We will also observe that 
Benedetti's exposition is mathematically sound,
quoting \textit{Euclid's Elements} where needed.

The following description  in modern words 
is based on the translation in Section~\ref{App~II}.

\begin{enumerate} 
\item[--]
Given are the perpendicular minor semiaxis~$ae$ and the major semiaxis~$ce$, as well as the radius~$oa$ of the larger circular arc, making the assumption that~$oa > oc$.
 \ \\[-4ex]
\item[--] 
 The objective is to construct an oval with axes as prespecified 
 above with~$o$ being a centre of the circular arc through~$a$.
  \ \\[-4ex]
\item[--]
Draw the perpendicular  to~$ao$ through~$o$
and let~$b$ be its (left) intersection point with the circle around~$o$ of radius~$ao$.\\
Draw the line~$bc$ up to the intersection point~$d$ with the same circle. (The intersection at~$d$ exists since we assumed~$oa > oc$.)\\
Draw the line~$do$ and let~$t$ be its intersection with~$ce$.\\
Since~$ob$ and~$ec$ are parallel, the angles~$tcd$ and~$obd$ are equal, as are the angles~$ctd$ and~$bod$.\\
By construction, the triangle~$bod$ is isosceles, so the angles~$dbo$ and~$bdo$ are equal;
hence also the angles~$tcd$ and~$cdt$ are equal, the triangle~$ctd$ is also isosceles,
and we have~$tc=td$ as required.
Benedetti's proof of this contains references to the corresponding propositions in \textit{Euclid's Elements}.\\ 
Draw the point~$x$ as the intersection of~$ce$ with its perpendicular 
through~$d$.\\
Draw the point~$h$ as the intersection of this perpendicular 
and the circle around centre~$t$ with radius~$td$.\\
Now the oval in the second quadrant is constituted by the arc~$\wideparen{cd}$
 with radius~$td$
and centre~$t$, and by the arc~$\wideparen{da}$   with radius~$oa$
and centre~$o$.
\\ 
Set the point~$z$ on the semiaxis~$ae$ such that~$ze=oe$.
\\
Finally, the centres of the oval are given by~$t$, $o$, $r$ and~$z$;
reflecting~$d$ and~$h$ through the line~$ae$ completes the construction of the oval.
\end{enumerate}



\captionsetup[subfloat]{labelformat=empty}

\subsection[Benedetti's oval~B2]{Benedetti's oval~B2\quad \normalfont{(Fig.\,\ref{Abb:Benedetti-B2})}}\label{Ssec:B2}

In the other part of the second letter, Benedetti constructs an oval given both semiaxes, as well as the length of the radii of the smaller circular arcs.~(Fig.~\ref{Abb:Benedetti-B2}) 
We refer to this construction as \textit{oval~B2}. 
Once again, the three parameters determine the oval uniquely. The translation in Section~\ref{App~II} forms the basis for the following description, which is similar to the previous construction and therefore shorter, again omitting Benedetti's citations of results from \textit{Euclid's Elements}.

\begin{enumerate}
\item[--]
Given are the minor (horizontal) axis~$qa$,
 the major (vertical) axis~$cp$  (satisfying~$qa<cp$),  meeting perpendicularly at~$e$,
and the circle of radius~$ct$ around~$t$ on the (upper) major semiaxis with $ct < ae$.
 \ \\[-4ex] 
\item[--]
 The objective is to construct an oval with axes and circular arc as specified above.
  \ \\[-4ex]
\item[--]
Draw the point~$x$ on~$ae$ so that~$xa=ct$.\\
Let~$n$ be the midpoint of the line~$tx$.\\
Let~$o$ be the intersection point of the perpendicular bisector of~$tx$ with~$qa$.\\
Then, $o$ is the required centre of the circular arc~$\wideparen{da}$, 
since $od=ot+ct=ox+xa=oa$,
whereas~$t$ is the centre of the circular arc~$\wideparen{cd}$. 
\end{enumerate}

Benedetti finishes with a remark: instead of constructing the perpendicular bisector of~$tx$ one could determine~$o$ as the point on the half-line~$eq$ for which the angle~$tox$ equals the angle~$txo$.
This also ensures the triangle~$otx$ is isosceles.


\subsection{Benedetti's encompassing circles}\label{Ssec:Bcircles}
%
%
In the first letter to Pizzamano, Benedetti solves the problem of constructing a
circle which encompasses two given circles.~\cite[pp.\,262/3]{Benedetti-1585}
The given circles may differ in size, and 
 three different cases of the positioning of the circles
(adjacent, overlapping, separate)  are considered.~(Fig.\,\ref{Fig:Benedetti-262-263-circles})
For these cases Benedetti suggests four solutions so that the 
encompassing circle touches both given circles.

The constructions of the encompassing circles for adjacent and overlapping circles~--
treated in Subsection~\ref{Sssec:B-262-oben}~(Fig.\,\ref{Abb:B-262-oben})
and in Subsection~\ref{Sssec:B-262-unten}~(Fig.\,\ref{Abb:B-262-unten})~-- 
 differ only in minor points.
The construction in Subsection~\ref{Sssec:B-263-oben}~(Fig.\,\ref{Abb:B-263-oben}) of a circle encompassing separate circles
yields exactly one circle, 
whilst a family of circles encompassing separate circles is derived in Subsection~\ref{Sssec:B-263-unten}~(Fig.\,\ref{Abb:B-263-unten}).
 
\begin{figure*}[t!]
\begin{center}
\subfloat[(a) \ adjacent circles]
 {
 \centering
   \includegraphics[angle=0,width=6.6cm]{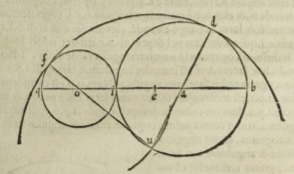}
\label{Abb:B-262-oben}
 \vspace{5mm}
 }
 \hspace*{6mm}
\subfloat[(b) \  overlapping circles]
 {
  \centering
   \includegraphics[angle=0,width=6.6cm]{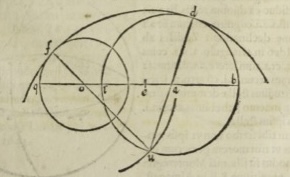}
\label{Abb:B-262-unten}
}
 \qquad\quad
\subfloat[(c) \ separate circles, special case]
 {
  \centering
   \includegraphics[angle=0,width=6.6cm]{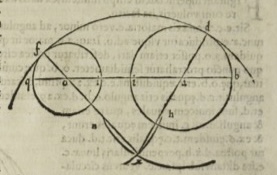}
\label{Abb:B-263-oben}
 }
  \hspace*{6mm}
\subfloat[(d) \  separate circles, general case]
 {
  \centering
   \includegraphics[angle=0,width=6.6cm]{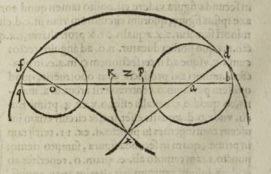}
\label{Abb:B-263-unten}
 }
 \caption{\small  Benedetti's four cases of given circles; details from~\cite[pp.\,262--263]{{Benedetti-1585}}.} 
 \label{Fig:Benedetti-262-263-circles}
 \end{center}
\end{figure*}

\subsubsection[Adjacent circles]{Adjacent circles\quad \normalfont{(Fig.\,\ref{Abb:B-262-oben})}}\label{Sssec:B-262-oben}

\begin{enumerate}
\item[--]
Given are on one line:
the radius~$qo$ of a (left) circle with centre~$o$,
the radius~$ab$ of a larger (right) circle with centre~$a$, and
the connection point~$i$ of both circles; $i$ is on the line through the centres~$o$ and~$b$ since the circles are assumed to touch each other there, for which Benedetti quotes Euclid~\cite[Prop.\,III.11]{Euclid-Clavius-1591-engl}.
 \ \\[-4ex] 
\item[--]
The objective is to construct an encompassing circle in the sense that
it touches the small circle at~$f$ and 
the large circle at~$d$.
 \ \\[-4ex] 
\item[--]
Let~$e$ be the point on the horizontal to the right of~$i$ such that~$qo = ie$.
Let~$u$ be an intersection point of
  the circle around~$o$ with radius~$eb$ and
    the circle around centre~$a$ with radius~$ab$.
 The existence of~$u$ follows since~$qo \leq ab$ implies that the sum of the radii is at least the distance of the centres,
$ ou + ua =  (ea + ab) + ab \geq  qo + ab= oa $.
\\
Now let~$f$ be the intersection point of the circle around~$o$
and radius~$qo$ and the  extension of the line segment between~$o$ and~$u$.
\\
Let~$d$ be the intersection point of the circle around~$a$ with
radius~$ab$ and the extension of the line segment between~$u$ and~$a$.
\\
By the same proposition of Euclid, the circle around~$u$ with radius
\[
fu = fo + ou = qo + eb = qo + ea + ab = ia + ab = ua + ab = ud
\]
touches both the smaller circle at~$f$ and the larger circle at~$d$.
\end{enumerate}

\captionsetup[subfloat]{labelformat=empty}

\subsubsection[Overlapping circles]{Overlapping circles\quad \normalfont{(Fig.\,\ref{Abb:B-262-unten})}}\label{Sssec:B-262-unten}

Here the case is considered where all conditions as in Subsection\,\ref{Sssec:B-262-oben}
remain unchanged, except that the given circles are no longer tangent but intersect (and~$i$ now denotes the intersection point of the larger circle with the line~$oa$).
The construction goes through without any changes and yields the desired result, 
so Benedetti just states ``The very same I say for intersecting circles''.

\subsubsection[Separate circles]{Separate circles\quad \normalfont{(Fig.\,\ref{Abb:B-263-oben})}}\label{Sssec:B-263-oben}

Before considering the next construction, observe that in Fig.\,\ref{Abb:B-263-oben} and in Benedetti's text, the two points on each given circle closest to 
each other have both been labelled with~$i$. For clarity, we will denote one of them by~$j$.

Now the case where the given circles are separate is studied.

\begin{enumerate}
  \item[--]
  Let two circles with centres~$a$ and~$o$ be given, and denote their intersection points with the line~$oa$ by~$b,i$ and~$q,j$, respectively,
  such that~$i,j$ are on the line segment~$oa$ and~$qo=oj$,~$ia=ab$.
   \ \\[-4ex] 
  \item[--]
  Assume~$2(ab + oq) \geq oa$
  which is the formal way to express Benedetti's informal restriction
  excluding that ``the distance between the two 
  given circles [be] so large that the resulting circles [in the next step] could neither touch each other nor intersect''.
\\  
  Let the point~$x$ be the intersection below the line segment 
  of  the circle around~$o$ with radius~$bi$
  and the circle around~$a$ with radius~$oq+ab$.
  Then we have\\
  \centerline{$xf = fo + ox = fo +bi = oq+ab+ai = xa+ad =xd$}
  and thus
    the circle around~$x$ with radius~$xf$
   yields the sought  circular arc which encompasses and touches the two given
   circles.
\end{enumerate}
 
\subsubsection[Separate circles revisited]{Separate circles revisited\quad \normalfont{(Fig.\,\ref{Abb:B-263-unten})}}\label{Sssec:B-263-unten}

In order to avoid the restrictive assumption in the previous subsection, Benedetti presents another construction for the case where the two given circles neither intersect or touch.
In fact, as Benedetti himself observes, this results in a general parametrized family of encompassing circles.

\begin{enumerate}
\item[--]
Let the given separate circles with centres~$o$ and~$a$ intersect the line~$oa$ outside the segment~$oa$ in points~$q$ and~$b$, respectively,
\\
and let~$z$ be the midpoint of the line segment~$qb$.
\\
Choose~$k$ and~$p$ on the line~$oa$ equidistant from~$z$ with 
$op>oz$, $ak>az$.
\\
Let~$x$ be an intersection point of
the circle around~$o$ with radius~$op$ and 
the circle around~$a$ with radius~$ak$.
($x$ exists since $op+ak > oz+az=oa$.) 
\\
Now, if $f$ denotes
 the intersection point of 
 the circle around~$o$ with radius~$qo$ 
and the extension of the line segment~$xo$, 
\\
and~$d$ denotes
 the intersection point of
 the circle around~$a$ with radius~$ab$
and the extension of the line segment~$xa$,
\\
then the circle around~$x$ with radius
\[
xf = qo + ox = qo + op = qp = bk = ak + ab = xa + ab = xd
\]  
touches both given circles
and the task is solved.
\end{enumerate}
 
Note that the choice of~$kz = zp$ parametrizes all circles encompassing the two given circles.

\subsection{The relationship between the two letters}
\label{Sec:relation}

Benedetti's results on encompassing circles are
 not about ovals.
In fact, as Benedetti writes in his letter, the problem of constructing such encompassing circles was posed by Pizzamano, and Benedetti presents him with a solution.
The problem and its solution belong to the field of classical Euclidean geometry, and it is not unusual for a well-educated nobleman like Pizzamano to be interested in such a topic.
Indeed, the question he poses is a much simpler version of the famous problem of Apollonius where a circle simultaneously touching three given circles is sought.

In the special case, however, where the two separate circles in Subsection~\ref{Sssec:B-263-unten} have equal radii, 
the corresponding construction yields the upper half of an oval -- which by
reflection gives a complete oval: the given circles form the smaller circular arcs, the encompassing circle one of the larger circular arcs.~(Fig.\,\ref{Abb:B-263-unten}, cf.~Fig.\,\ref{Abb:Clavius-1591-150})
We thus obtain a construction of an oval where the radii of the oval's smaller circular arcs are given as~$qo=ab$, as well as their centres~$o$ and~$a$; equivalently, given the major semiaxis~$qb$ together with those radii the smaller arcs' centres are specified.
The choice of the distance~$kz=zp$ then  determines the radius of the encompassing circle, $xf = qo+oz+zp$ with~$2oz = oa$.
As this is the radius of the oval's larger circular arcs, we may as well assume that to be given.
Thus, the construction in Subsection~\ref{Sssec:B-263-unten} yields the family of all 4-centre ovals for given radii of the circular arcs and the major semiaxis.
We do not know of any other construction of an oval given these parameters, cf. e.g.\,\cite{Mazzotti-2019}.

\begin{figure}[t!]
\centering
\includegraphics[width=0.95\columnwidth]{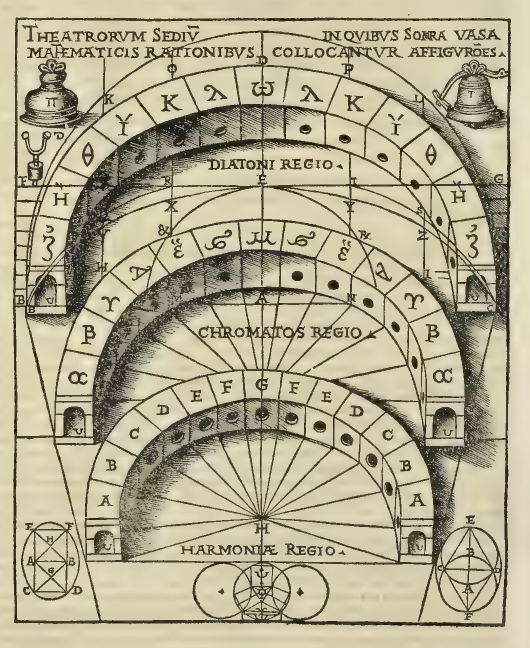}
\\[2ex]
\caption{\small 
Detail from Cesariano's edition of Vitruv from~1521~\cite[p.\,LXXX verso]{Vitruv-1521-Cesariano-engl}.\\[3ex]
} 
\label{Abb:Cesariano-LXXX}
\end{figure}

Benedetti himself does not mention this connection between the construction of encompassing circles and the construction of ovals -- not least because apparently he was not familiar with the notion of an oval.
But the fact that the two letters were both addressed to Pizzamano, and were included in the \textit{Diversarum speculationum} with the one about ovals immediately following the one concerning encompassing circles, suggests that Benedetti first answered Pizzamano's question about the construction of encompassing circles; and he then observed that for given circles of equal radii the last construction (after reflection) results in ``a figure which superficially resembles an ellipse'' (which we now call an oval).
In fact, the letter on ovals does not include any hint that it directly answers a letter by Pizzamano, but rather reads as if Benedetti simply wanted to tell Pizzamano about his curious observation.

The oval constructions in the second letter differ from the general construction of encompassing circles in an important aspect, though: the former assume that both semiaxes and only one of the circular arcs' radii is given.
It appears that Benedetti, viewing ovals as figures resembling ellipses which are determined by their semiaxes, also wanted to directly construct ovals with prescribed semiaxes.
Nonetheless, the resemblance (after rotation) of the construction in Subsection~\ref{Sssec:B-263-unten}~(Fig.\,\ref{Abb:B-263-unten}) and the construction of oval~B2 (Fig.~\ref{Abb:Benedetti-B2}) is striking.
It stands to reason that Benedetti obtained the latter construction of oval~B2 starting from the earlier one.

In fact, in the construction of~Fig.\,\ref{Abb:B-263-unten} the larger radius~$xf$ is composed of the smaller (given) circle's radius~$qo$ and the distance~$oa$ between the circles' centres.
Therefore, in the situation of oval~B2~(Fig.\,\ref{Abb:Benedetti-B2}), it suggests itself to mark the smaller radius~$ct$ off the minor semiaxis~$ea$, resulting in~$x$; now, the larger circle's centre~$o$ needs to be at the same distance from~$x$ as from the smaller circle's centre~$t$, which Benedetti realises by intersection with the bisector of the line segment~$tx$.
Thereby one obtains Benedetti's construction of oval~B2 in Subsection~\ref{Ssec:B2} from the one in Subsection~\ref{Sssec:B-263-unten}.

The construction of oval~B1 (Fig.~\ref{Abb:Benedetti-B1}) appears to follow a rather different line of thought.



\section{Benedetti's ovals in the literature}\label{Sec:Reception}
In the following, we aim to create a chronological list of works treating Benedetti's ovals. However, for the sake of completeness, and because later authors often refer to Serlio's ovals, we make an exception by presenting the precursors to Serlio's ovals in the first subsection.

Recall that Benedetti sent his results on ovals in a letter to
Pietro Pizzamano who died in~1571, so the letter must have been written at least~14 years before the first edition of the \textit{Diversarum Speculationum} appeared in~1585.

In our  review of the impact Benedetti's work on ovals had,
we attempt to determine how closely the authors adhere to the constructions of oval~B1 and oval~B2, and whether their presentations are correct. As we will see, this is not always the case.

\subsection{Precursors to Serlio's ovals}

In 1521, the Italian painter and architect Cesare Cesariano (1483–1543) translated and commented on Vitruvius' \textit{De Architectura Libri Decem} into Italian~\cite{Vitruv-1521-Cesariano-engl}. In one of his figures (Fig.~\ref{Abb:Cesariano-LXXX}), one finds, almost hidden in the lower left and right corners, two ovals which were later presented by Serlio, namely ovals~S3 and~S4.  These ovals were published~24~years before Serlio's treatise was
 printed in~1545.
 
Between~1532 and~1536, the Italian architect and painter
Baldassare Peruzzi (1481--1536)   
presented some sketches of two ovals 
which again coincide with ovals~S3 and~S4. (See Uffizi archive~531,
also  shown in~\cite[p.\,12, Abb.\,1]{Lotz-1955-engl}.)
These ovals were sketched  at least~9~years before Serlio's treatise was
 printed.
 Note that Serlio worked under Peruzzi.~\cite{Duvernoy-2015-engl}

%
\subsection{Christophorus Clavius 1591 and 1604} \label{Ssec:Clavius}
Christophorus Clavius~(1538--1612) was a mathematician and Jesuit priest at the Collegio Romano. 
He was probably the first to comment on the significance of Benedetti’s work on ovals.
In his \textit{Geometria Practica},
 which first appeared in~1604 (and in subsequent editions), 
ovals are treated in Liber Octavus, Problema~33,   Propositio~47.~
\cite[pp.\,421--424]{Clavius-1604-engl}

\begin{figure*}[t!]
   \begin{center}
 \subfloat[(a) \ Clavius' generalization of oval~S1]  
 { 
  \centering
   \includegraphics[angle=0,width=4.2cm]{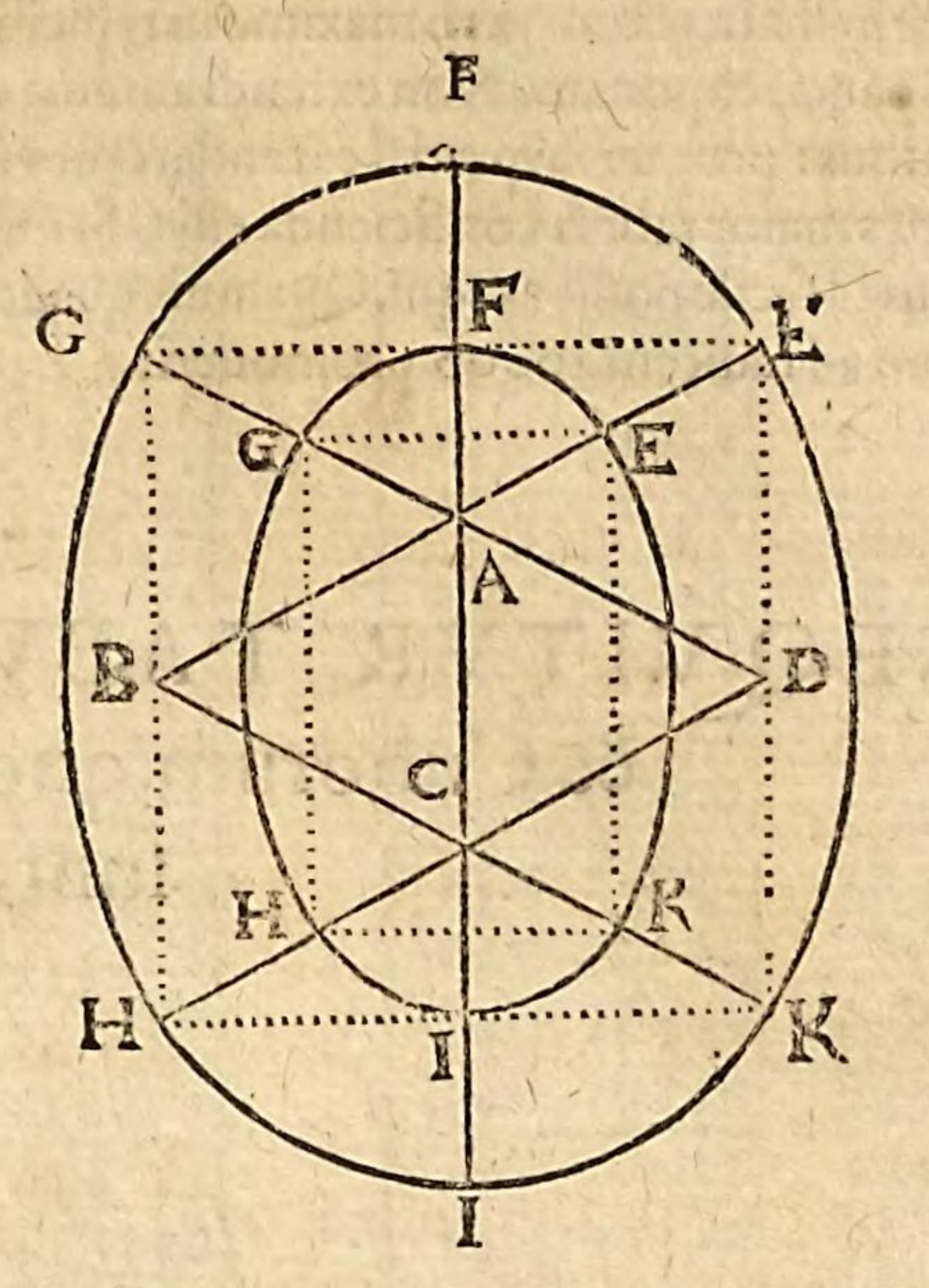}
\label{Abb:Clavius-1611-Geom-pract-S219o}
 } 
\hspace*{2ex}
 \subfloat[(b) \  Clavius' version of oval~B1]  
 { 
  \centering
   \includegraphics[angle=0,width=4.2cm]{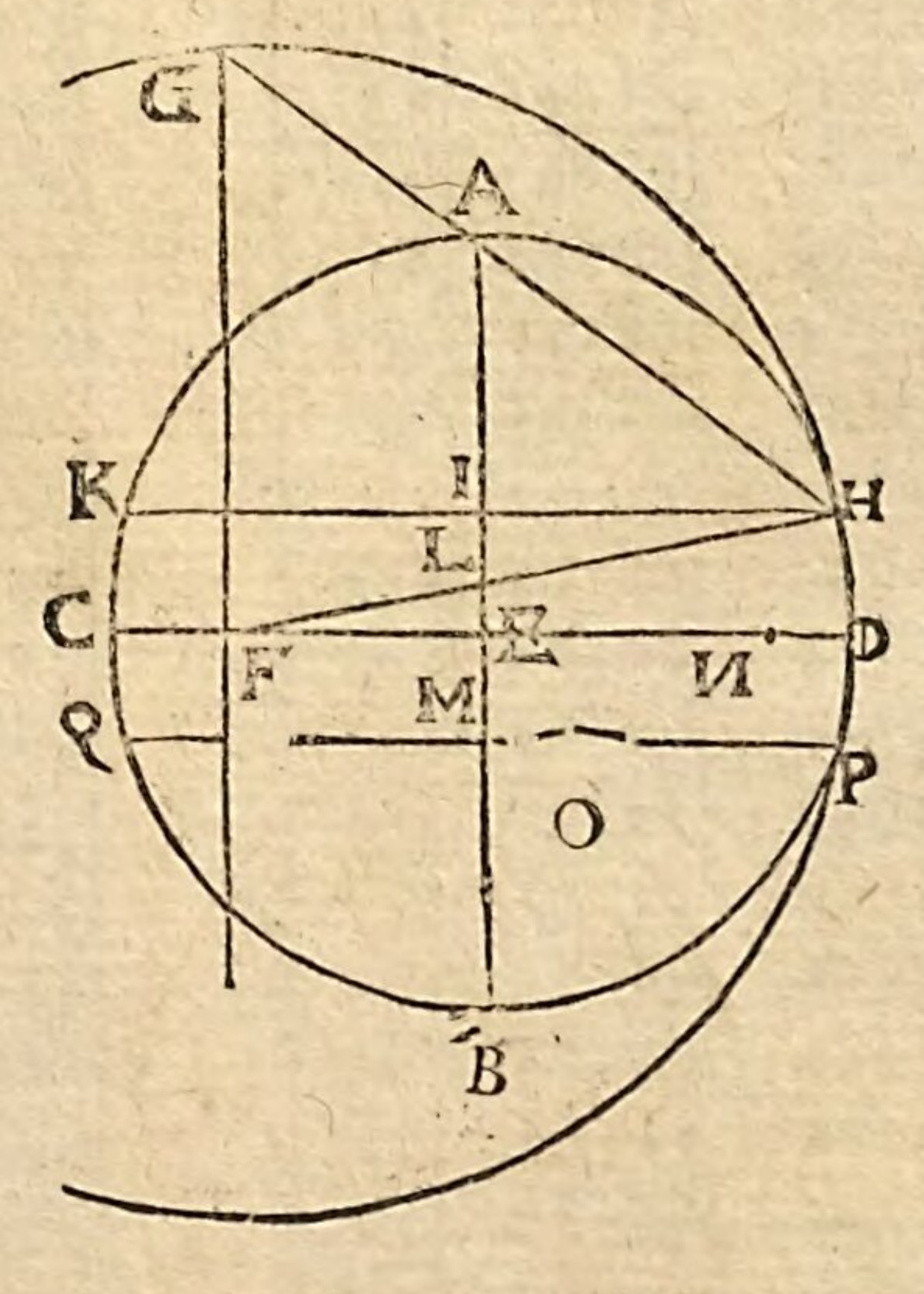}
\label{Abb:Clavius-1611-Geom-pract-S219u}
 } 
 \hspace*{2ex}
 \subfloat[(c) \  Clavius' version of oval~B2]
 {
 \centering
   \includegraphics[angle=0,width=4.65cm]{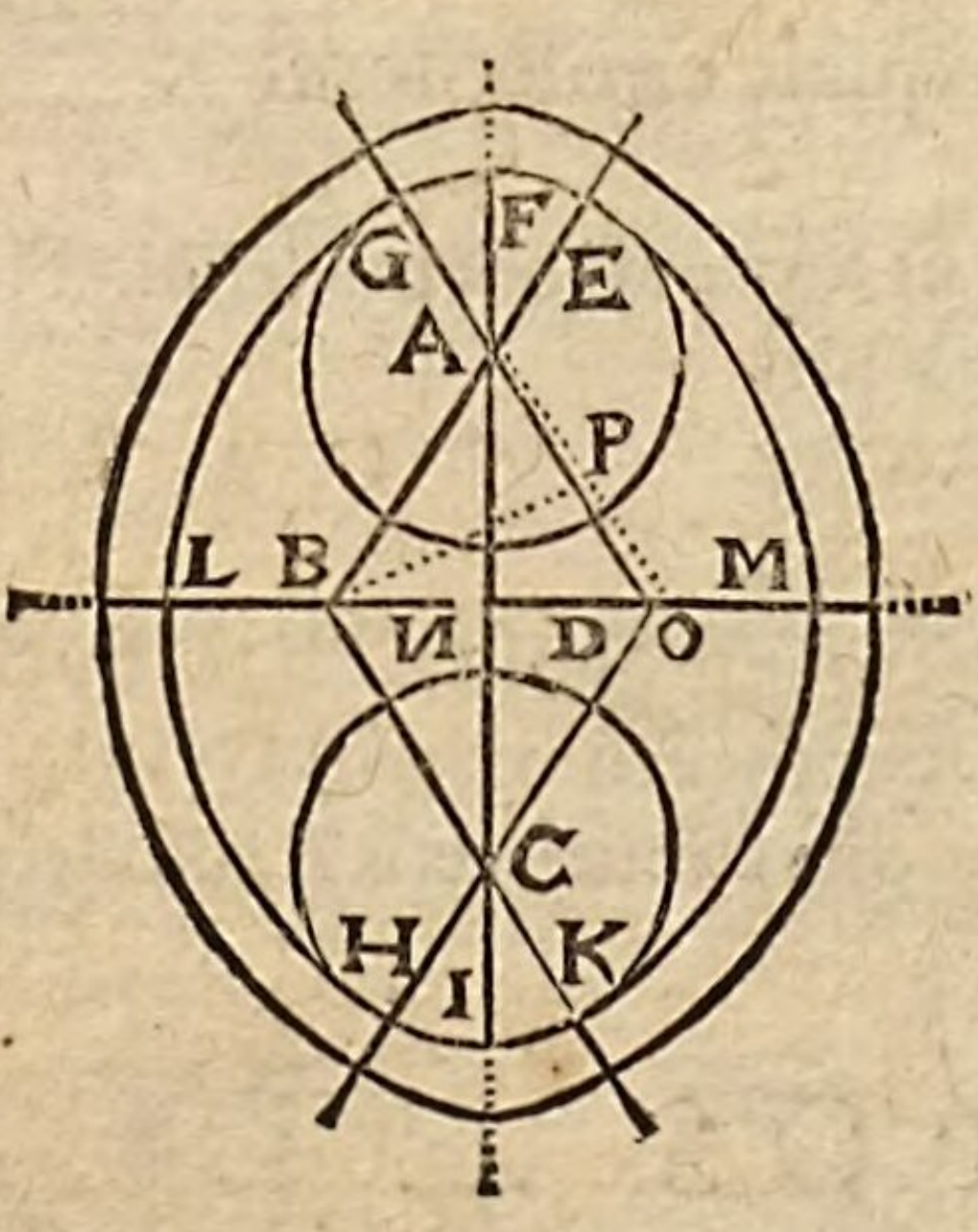}
\label{Abb:Clavius-1611-Geom-pract-S218}
 }
 \vspace{1.1ex}
 \caption{\small  Details from~\cite[pp.\,218--219]{Clavius-1611-engl} as in Clavius 1604~\cite[pp.\,421--423]{Clavius-1604-engl}.}
 \label{Fig:Clavius}
 \end{center}
\end{figure*}

Clavius starts by studying variants of Serlio's oval~S1 (without giving credit), 
but going beyond Serlio's construction, 
he no longer requires the starting triangle~$ABC$ to be equila\-teral.~(Fig.\,\ref{Abb:Clavius-1611-Geom-pract-S219o})
From an art theoretic viewpoint, this represents an emancipation of form.
Serlio’s ovals are constituted by ideal, regular forms such as circles, squares, equilateral triangles, or the vesica piscis.
  Regularity was seen as necessary for beauty, a concept that traces back to Aristotle\footnote{Aristotle:
``A science is most precise if it abstracts from movement, but if it takes account of movement, it is most precise if it deals with the primary movement, for this is the simplest; and of this again uniform movement is the simplest form.''~\cite[1078a]{Aristotle-1927}
}
 and remained valid during Baroque times. 
 From a purely geometric viewpoint, however, isosceles triangles are less restrictive than equilateral ones, allowing for greater flexibility in design, making it possible to obtain any 4-centre oval by this construction.
 But Clavius also notes that more ``pleasing (\textit{venustiores})'' figures are obtained when using an isosceles triangle, ``as artists use to do'', and refers to the two ovals in Fig.~\ref{Abb:Clavius-1611-Geom-pract-S219o} where that is the case.

%
%

\captionsetup[subfloat]{labelformat=empty}

However, Clavius also states that 
``in this manner, no figure with a given length and width could be drawn''.~\cite[p.\,423]{Clavius-1604-engl}
Continuing, he quotes Benedetti for having obtained solutions to this problem and constructs 
oval~B1. 
Indeed, when rotated anti-clockwise by~90°, Fig.\,\ref{Abb:Clavius-1611-Geom-pract-S219u}
 coincides with
 Fig.\,\ref{Abb:Benedetti-B1} 
 with~$F$ as the centre of the larger circular arc,
 $L$ the centre of the smaller one,
 and~$AE,ED$ as the given semiaxes.

Clavius also constructs oval~B2~\cite[p.\,424]{Clavius-1604-engl} rigorously,
compare
Fig.\,\ref{Abb:Clavius-1611-Geom-pract-S218}
with
Fig.\,\ref{Abb:Benedetti-B2},
though (deviating from Benedetti) falsely claiming the first construction not to be general, thus demanding another one.

\begin{figure*}[t!]
\begin{center}
\subfloat[(a) \  Clavius~1591]
 {
 \centering
   \includegraphics[angle=0,width=6.6cm]{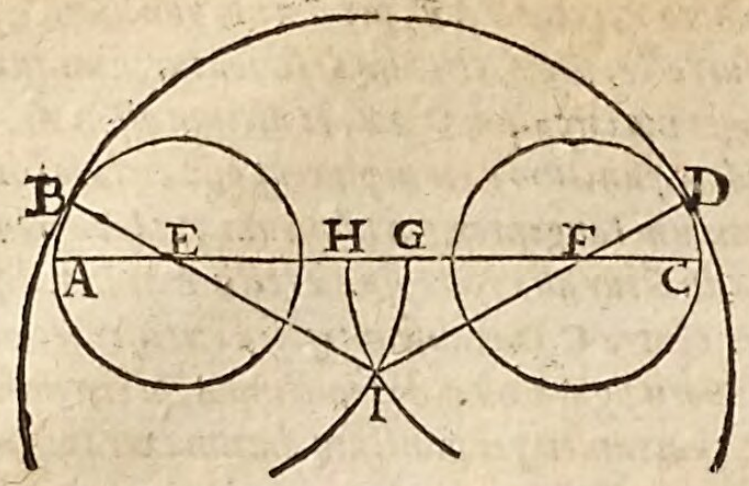}
\label{Abb:Clavius-1591-150}
 \vspace{5mm}
 }
 \hspace*{6mm}
\subfloat[(b) \  Ardüser~1627]
 {
  \centering
   \includegraphics[angle=0,width=6.6cm]{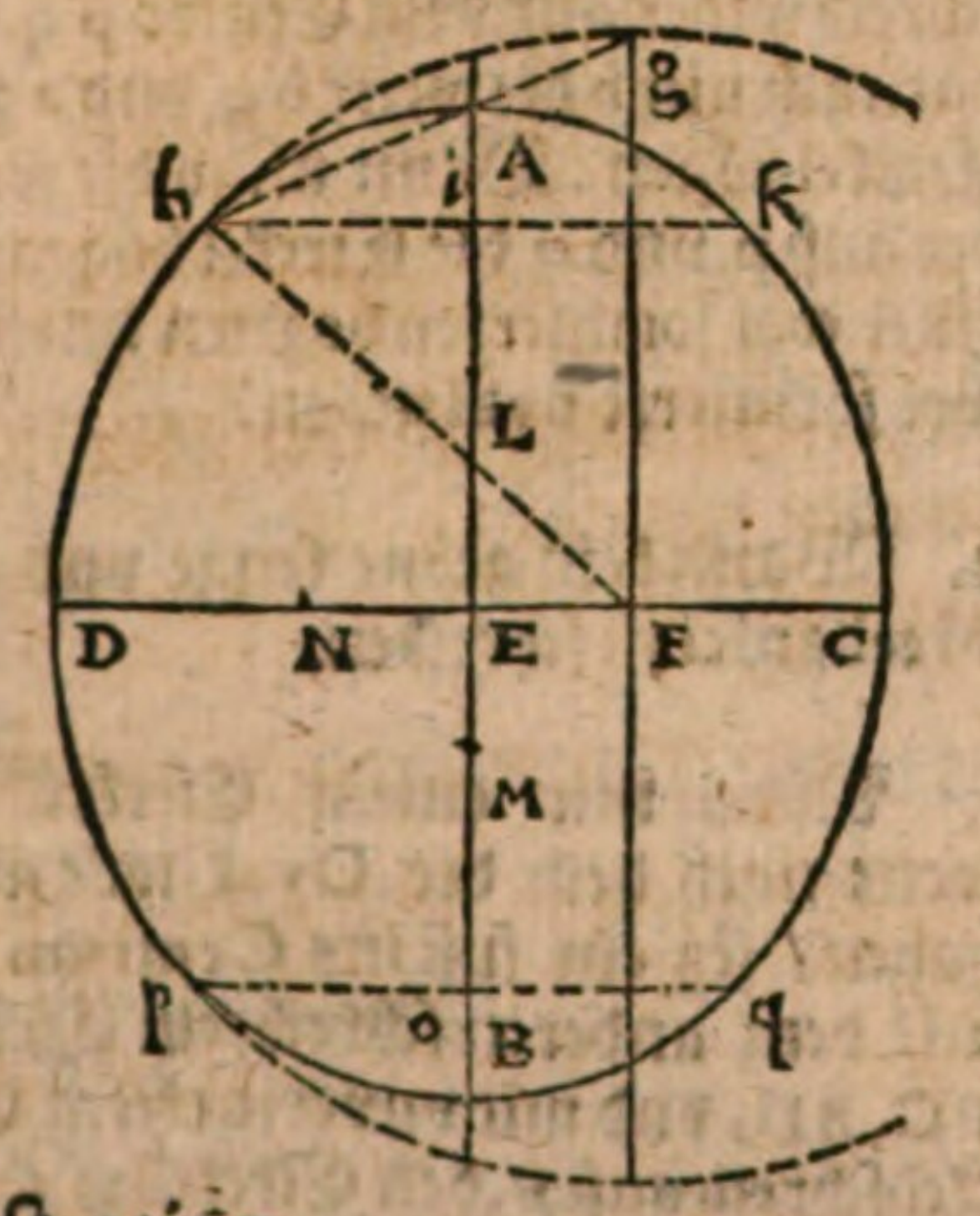}
\label{Abb:Ardueser-1627-S134}
}
 \caption{\small  Details from Clavius~\cite[p.\,150]{Euclid-Clavius-1591-engl} and Ardüser~\cite[p.\,134]{Ardueser-1627}}
 \label{Abb:Clavius-1691-Ardueser}
\end{center}
 \end{figure*}

In his famous Euclid edition from~1591, Clavius treats Benedetti's construction of a circle encompassing two given circles, too.~\cite[pp.\,149--150]{Euclid-Clavius-1591-engl}
There, he follows Benedetti (see Subsection~\ref{Ssec:Bcircles})
in his text and illustrations  which are almost identical to Benedetti's;
he not only quotes Benedetti, but furthermore discusses the case where the two given circles are of identical size, i.e., where one almost obtains oval~B2, cf. Subsection~\ref{Sec:relation} and Fig.\,\ref{Abb:Clavius-1591-150}.
However, presenting the contents of the two letters in different works, Clavius does not comment on the relationship of these constructions.

\captionsetup[subfloat]{labelformat=empty}

\subsection{Matthäus Beger 1623}

The only known translation of one of Benedetti's letters to Pizzamano into a modern language, namely into German, has been manufactured by Matthäus Beger in~1623.
He lived at the free imperial city Reutlingen where he was born in~1588, and died in~1661. During his apprenticeship as a cloth shearer at Ulm he came into contact with the mathematical practitioner and engineer Johann Faulhaber who introduced him to the mathematical sciences.
In his spare time, Beger then studied mathematics as well as the languages this required by himself,  in particular Latin.
Supported by his business, he assumed several offices for the city: in the militia, as negotiator during the Thirty Years' War, as treasurer and finally as the city's mayor.
He wrote more than thirty manuscripts in folio on mathematical topics, mainly containing translations and excerpts from different authors for his fellow citizens to read and learn from.
In~1652 he donated these volumes along with several hundred prints in his possession to the city, thereby founding Reutlingen's city library where 
most of them remain.~\cite{Brinkhus-2022,Sommer-1978}

Beger's manuscript~21 (Friderich catalogue 1881), entitled \textit{Secreta Mathematica} contains translations of nine works of different authors, number 6 (pp.\,244r--438r) being Benedetti's \textit{Diversarum speculationum}, in the volume's table of contents called \textit{Speculationes Mathematicae}, as Beger used the 1599 edition named \textit{Speculationum Liber} which he had borrowed from Matthias Bernegger, professor at the University of Strasbourg.
The translation is dated from~1623. It comprises a large part of the \textit{Diversarum speculationum}; Beger notes (page~224v) that he deemed it unnecessary though to include every of its propositions as much of the material is dealt with better by other authors.
Nonetheless, out of Benedetti's~118 letters he translated~39, in particular as number~7 (pp.\,399r--399v) the one to Pizzamano on ovals~B1 and~B2 ``Daß ein Superficialische Figur einer Ellipsi gleichförmig, auß denn gegebnen Axibus mit einem Circkel delineirt werden mag''.
Unfortunately, we do not know whether anyone ever took note of Beger's translation of this letter, as it was never printed.
 
\subsection{Johann Ardüser 1627}    
In 1627, Johann Ardüser (1585--1665), a Swiss mathematician and  fortress engineer,
published his \textit{Geometriae, Theoricae et Practicae} (in German).
In Book 5, Section~IIII~\cite[pp.\,133v--135r]{Ardueser-1627} he discusses ovals (which he identifies with ellipses, though).
He begins with the construction of Benedetti's oval B1 together with a short proof.
Indeed, if you rotate his figure depicted in Fig.\,\ref{Abb:Ardueser-1627-S134} counterclockwise by 90° and then flip it along the horizontal axis, it becomes oval~B1, as shown in Fig.\,\ref{Abb:Benedetti-B1}.
Afterwards, he shows an approximative method to construct a similarly shaped figure; finally, he presents ovals S4, S2, and S1.
We note that Ardüser cites neither Serlio nor Benedetti but Clavius in his preface, so he may have taken these oval constructions from Clavius' \textit{Geometria Practica}.
\\

\captionsetup[subfloat]{labelformat=empty}
 \begin{figure*}[t!]
\begin{center}
\subfloat[(a)  \  Bettini 1645]
 {
 \centering
   \includegraphics[angle=0,width=6.6cm]{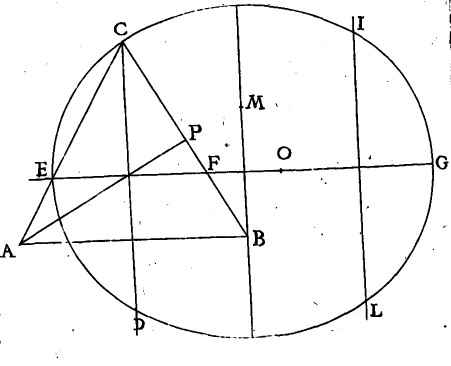}
\label{Abb:Bettini_1645_Apiaria-106}
 \vspace{5mm}
 }
 \hspace*{6mm}
\subfloat[(b)   Harsdörffer~1653]
 {
  \centering
   \includegraphics[angle=0,width=6.6cm]{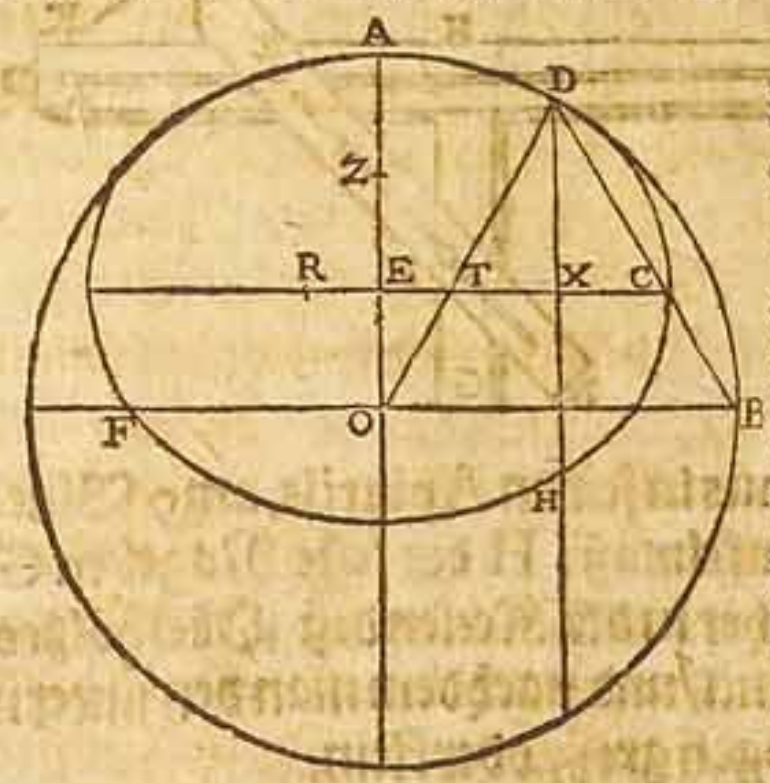}
\label{Abb:Harsdoerffer3-1653-S180}
}
 \caption{\small  Details from Bettini~1645~\cite[Tl.\,1, Apiar.\,III, p.\,106]{Bettini-1645-55-engl}
 and Harsdörffer~1653~\cite[p.\,180]{Harsdoerffer-1653}
 }
 \label{Abb:Bettinie-Harsdoerffer}
\end{center}
 \end{figure*} 

\subsection{Mario Bettini 1645}
Mario  Bettini~(1583--1657) was a  Jesuit priest and professor of mathematics at Parma who published his \textit{Apiaria Universae  Philosophiae Mathematicae} (1642--1655) where he discusses ovals  in the section titled
\textit{Pseudoellipsis, sive oviformis figurae descriptio compendio perfecta, \& theorice explicata}.~\cite[vol.\,1, sect.\,III.X, prop.\,XIII, pp.\,106/7]{Bettini-1645-55-engl} 
Prima facie,  one might think that his oval depicted in
Fig.\,\ref{Abb:Bettini_1645_Apiaria-106}
resembles   Benedetti's oval~B1 in
Fig.\,\ref{Abb:Benedetti-B1}, 
as both constructions invoke an isosceles triangle.
However, Bettini does not presume the semiaxes to be given.
Instead, he starts with given length~$AB$ (which will later be the radius
of the larger circular arc) and draws an equilateral triangle~$ABC$ (thus losing generality).
The centres of the oval are then determined via the centre of~$ABC$ 
as~$F$,~$B$ and their reflections~$O$, and~$M$. 
One may wonder whether this was obtained through an attempt to understand the construction of oval~B1 by merely looking at it.

Furthermore,  Bettini suggests another construction of an oval and 
claims that ``the entire construction of this figure can, however, also be done more simply and with less effort''.~\cite[p.\,107]{Bettini-1645-55-engl}
Here, Bettini starts with the longer semiaxis, and the centres of the arcs as given which is essentially equivalent to 
the data given for 
Clavius' generalisation of Serlio's oval~S1; cf.~Subs.\,\ref{Ssec:Clavius}.
\\


\subsection{Georg Philipp Harsdörffer 1651}
Georg Philipp Harsdörffer~(1607--1658) was a German poet and polymath who published more than 50 books. His two volumes, \textit{Der Mathematischen und Philosophischen Erquickstunden Zweyter Theil}~\cite{Harsdoerffer-1651} and \textit{Dritter Theil}~\cite{Harsdoerffer-1653}, published in~1651 and~1653, became very popular and were often reissued. An English translation of the title might read \textit{The Mathematical and Philosophical Hours of Delight, Second Part} and \textit{Third Part}, which is meant to be entertaining.
\footnote{Note that the first part was written by Daniel Schwenter, professor of mathematics at Altdorf, and published posthumously in~1636.}



When he treats oval~B1~\cite[p.\,180]{Harsdoerffer-1653}, Harsdörffer quotes Benedetti and provides a figure (Fig.\,\ref{Abb:Harsdoerffer3-1653-S180}) that is, up to reflection, Benedetti's Fig.\,\ref{Abb:Benedetti-B1}.
However, his description is not convincing. He does not mention that~$O$ must fulfill~$AO > OC$, as discussed in Subsection~\ref{Ssec:B1}. Instead, he fixes~$O$ by setting~$EO = AE/2$ which violates the condition if~$OC$ is too large. Moreover, he seems not to have understood that~$O$ can be chosen freely (if~$AO$ is large enough). This free choice of~$O$ is important for architectural designs. Also, he does not mention the core feature of the construction, which is that the triangle~$TCD$ is equilateral. Therefore, his description does not qualify as a 
rule for constructing oval~B1.

\subsection{Abraham Bosse 1665}
%
Abraham Bosse~(ca.\,1604--1676)
was a well respected engraver und etcher of the~17th century.
In his \textit{Trait{\'{e}} des pratiques g{\'{e}}ometrales et perspectives}
from~1665, he presents~-- without citing Benedetti~--
 a construction of Benedetti's oval~B2. 
 %
To see that  Bosse's figure depicted in~Fig.\,\ref{Abb:Bosse-Oval_1665-T14-fig3} 
coincides 
with Benedetti's  figure in Fig.\,\ref{Abb:Benedetti-B2},
rotate the former  by~90° clockwise and identify 
the nomenclature as~$E\,\widehat{=}\, c$, 
  $I\,\widehat{=}\, t$, 
      $O\,\widehat{=}\, x$, 
        $S\,\widehat{=}\, n$, 
and the others accordingly.
Then both figures agree.


Bosse assumes that the major semiaxis~$EL$,
 the minor semiaxis~$ML$,
 and the radius~$EI$  of the small circular arcs are 
 arbitrarily given under the assumption~$EI<ML$.
 In passing we correct that~$N$ is not the end of the minor axis,
 as Bosse writes; but this does not affect Bosse's construction. 

The most important mathematical feature of the construction,
namely that the triangle~$INO$ is isosceles
which implies~$GN=MN$, is not mentioned explicitly. 
It could be proved as Benedetti did by referring to Euclid's \textit{Elements};
however, in the~17th century this might have been `common mathematical knowledge', not worth quoting. 

Compared to Benedetti, 
Bosse adds a discussion about the situation when~$MO$ becomes larger 
so that~$O$ moves to~$L$. 
In this case, the angle~$NOI$ approaches~90° and the distance~$LN$
tends to~$\infty$. Bosse writes:
``car si le point~$O$ est fort proche de~$L$, ayant tiré~$O\,I$, 
l'angle~$NOI$ fera presque droit, \&~la ligne~$S\,N$ rencontrera~$M\,N$
fort loin aus dessous de~$L$''.
And he ends with 
``ce qui seroit incommode''---``this would be inconvenient''. 
Indeed, a large~$LN$  means a flat arc~$GMH$,
and this might not be wanted because it does not yield a `proper' oval.


\subsection{Tomás Vicente Tosca 1712}
\begin{figure*}[t!]
\begin{center}
\subfloat[(a) \ Bosse 1665]
 {
 \centering
   \includegraphics[angle=0,width=7.2cm]{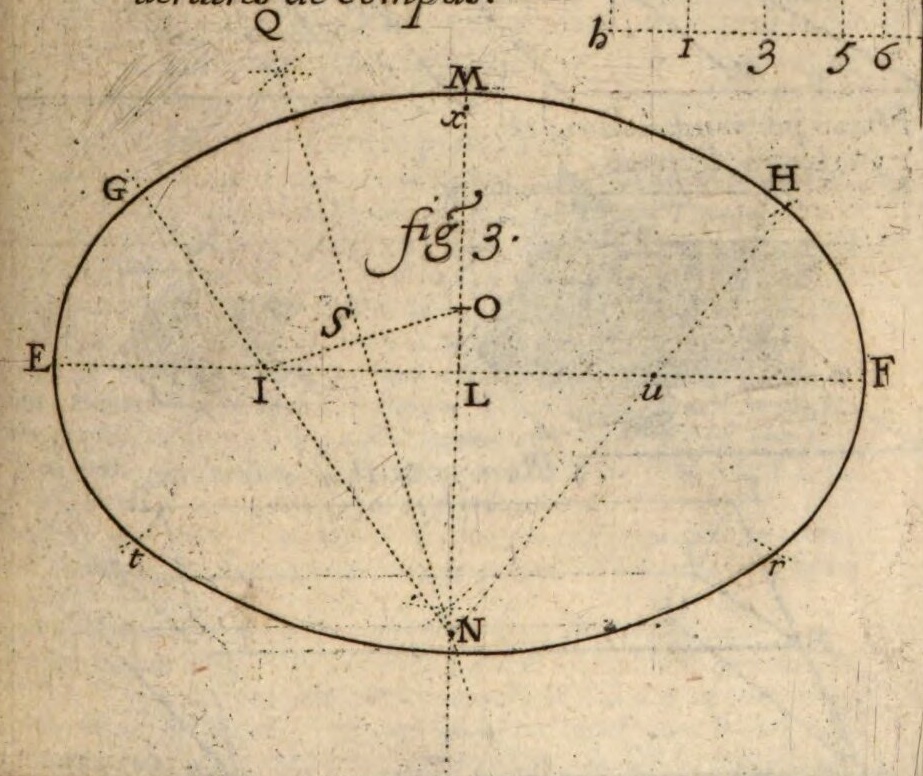}
\label{Abb:Bosse-Oval_1665-T14-fig3}
 \vspace{5mm}
 }
 \hspace*{4mm}
\subfloat[(b) \  Tosca~1712]
 {
  \centering
   \includegraphics[angle=0,width=6.3cm]{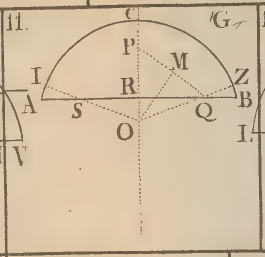}
\label{Abb:Tosca-1794-S129-Fig11}
}
 \caption{\small  Details from Bosse~\cite[p.\,64 \&\,Plate\,14]{Bosse-1665}
 and   Tosca~\cite[Folio~108, Fig.\,11]{Tosca-1727}.
 }
\end{center}
 \end{figure*}

The Spanish architect, priest and mathematician Tomás Vicente 
Tosca~(1651--1723)
published his  \textit{Compendio Mathematico, Vol.\,5} 
in~1712.~\cite{Tosca-1712}
This volume includes the oval depicted in  Fig.\,\ref{Abb:Tosca-1794-S129-Fig11}.
It is this construction which Ana L{\'{o}}pez Mozo claims to be ``possibly [\ldots] the first, general construction for drawing an oval for any given proportion''~\cite[p.\,580]{Lopez-Mozo-2011}; cf.\,Section~\ref{Sec:intro}.
 If this figure is turned 90° counterclockwise and then reflected through the vertical axis, one arrives at Benedetti's oval~B2,
  as shown in Fig.\,\ref{Abb:Benedetti-B2} with the associated nomenclature:
$M\,\widehat{=}\, n$, 
  $P\,\widehat{=}\, x$, 
          $Q\,\widehat{=}\, t$, 
      $O\,\widehat{=}\, o$, 
and the others accordingly.
Tosca describes his figure as
\begin{quote}
  ``\textit{Mode 2.}
  Let $AB$ be the horizontal diameter of the arch (fig.\,11), and $C$ the height it ought to have.
  Cut off arbitrary but equal pieces $AS$, $BQ$, $CP$. Draw the line $PQ$ which is divided in two equal pieces by the perpendicular $MO$, the latter intersecting the prolonged $CR$ in $O$.
  From $O$ through the point $Q$ draw the line $OQZ$, and $OSI$ through $S$; and from $O$ with the distance $OC$ construct the arch $ICZ$, and from $Q$ with distance $QZ$ construct the arch $ZB$, and from $S$ the arch $IA$.
  Thus, the interior side of the arch will be formed, and from the same centres the outer side will be produced.''~\cite[p.\,101]{Tosca-1712} 
\end{quote}
This description is complete except for Tosca's omission of the assumption~$2\,CR<AB$ (that is~$qa<cp$ in Benedetti's nomenclature) and of the proof of correctness.

\captionsetup[subfloat]{labelformat=empty}

\subsection{Giovanni Bordiga 1926}
The mathematician Giovanni Bordiga (1854--1933) was professor at the University of Padua and president of the \textit{Ateneo Veneto di Scienze, Lettere ed Arti}.
His biography of Benedetti~\cite{Bordiga-1926} is well-cited (cf. the references in Section~\ref{sec:benedetti}); in particular it appears to be one of the few commenting on the letters to Pizzamano.
Unfortunately, he writes on the first letter \cite[IV.10]{Bordiga-1926}:
``This is an unimportant problem, too. Indeed, the author shows how a circle can be constructed that is tangent to two other given ones;
but he does not say that the problem is underdetermined, neither does he look for the locus of the varying circle's center which satisfies the given conditions.''
This statement is in stark contrast with Clavius' more contemporary adoption of this letter's contents, see Subsection~\ref{Ssec:Clavius}.
Had Bordega carefully read the letter's last part, he might have noticed that Benedetti was indeed aware of the fact that the problem was underdetermined, and that he provided a general solution of varying circles satisfying the conditions.

The letter on ovals is treated with similar disparagement \cite[IV.11]{Bordiga-1926}: ``Nor should the following letters sent to him [...] be noted.''
These early pejorative comments by an authority on Benedetti may be one of the reasons why Benedetti's results on ovals have not yet been given much consideration in the literature on the history of mathematics and architecture.

\subsection{Angelo~Alessandro Mazzotti 2017}
\begin{figure*}[t!]
\begin{center}
\subfloat[(a) \  Mazzotti 2017: ``Construction~1a~-- Bosse's method'']
 {
 \centering
   \includegraphics[angle=0,width=6.6cm]{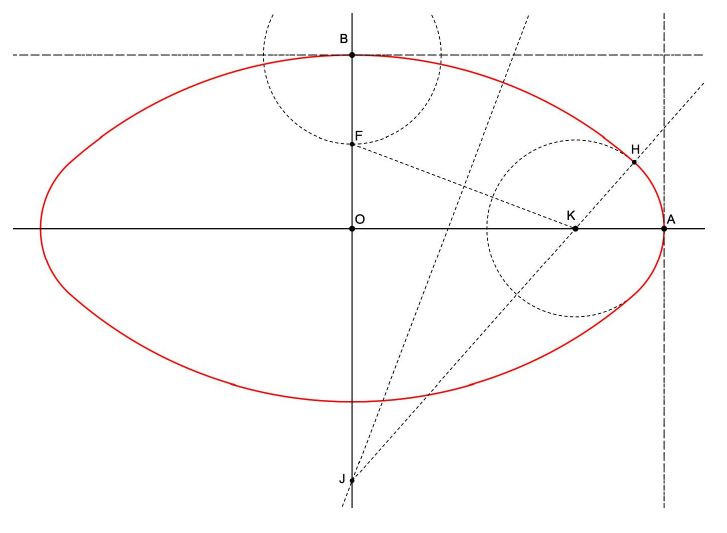}
\label{Abb:Mazzotti-2019-2nd-Fig31}
 }
 \hspace*{5mm}
\subfloat[(b) \  Mazzotti 2017: ``Construction~3a'']
 {
  \centering
   \includegraphics[angle=0,width=6.6cm]{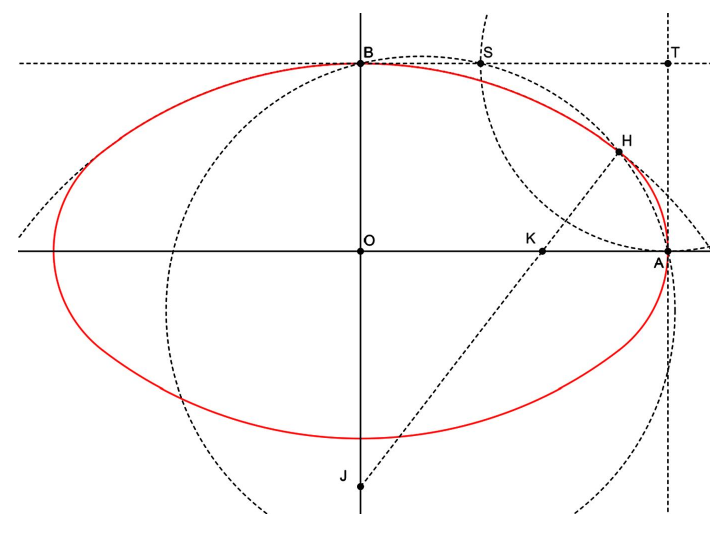}
\label{Abb:Mazzotti-2019-2nd-Fig31-3a}
}
\caption{\small 
Mazzotti's figures, which in the first edition of~2017 coincide with the ones in 
the second edition of~2019~\cite[p.\,20, Fig\,3.1 and p.\,25, Fig.\,3.5]{Mazzotti-2019}.
\label{Abb:Mazzotti-2019-2nd-Fig31u35}
}
\end{center}
 \end{figure*}

Angelo Mazzotti considers what he calls
 ``Bosse's construction'' twice.
First, in his \textit{Construction~1a} (see Fig.\,\ref{Abb:Mazzotti-2019-2nd-Fig31})
he assumes, as Benedetti does for oval~B2, and as Bosse does,
 that  the length~$O\!A$ of the major semiaxis,
 the length~$O\!B$ of the minor semiaxis,
 and the radius~$O\!A-O\!K$  of the small circular arcs are 
 arbitrarily given, under the assumption~$0<O\!A-O\!B<O\!K<O\!A$.
 The second inequality is, in Bosse's notation, equivalent to~$M\!O<M\!L$,
which is exactly what Bosse assumes.
However, Mazzotti does not explain where his inequality assumption is needed.

 If Mazzotti's 
 Fig.\,\ref{Abb:Mazzotti-2019-2nd-Fig31}
 is flipped along the vertical axis,
  it coincides with Bosse's Fig.\,\ref{Abb:Bosse-Oval_1665-T14-fig3}. 
In summary, Mazzotti refers to Bosse's construction where he is, in fact, reproducing Benedetti's construction of oval~B2.

Secondly,  
 Mazzotti presents another \textit{Construction~3a}
  (see Fig.\,\ref{Abb:Mazzotti-2019-2nd-Fig31-3a})
which he also refers to as ``Bosse's construction'',
and assumes, as Benedetti does for oval~B1,
 that  the length~$O\!A$ of the major axis,
 the length~$O\!B$ of the minor axis,
 and the radius~$BJ$  of the large circular arcs are 
 arbitrarily given under the assumption~$O\!J>  (O\!A^2-O\!B^2)/2O\!B$.
At first glance, one might think that the construction of oval~B1 follows, but in fact, it is different: 
 Mazzotti presents  a variant of the construction of oval~B2
 based on a very similar reasoning.
It can easily be verified that Mazzotti's inequality~$O\!J>  (O\!A^2-O\!B^2)/2O\!B$
 is equivalent to Benedetti's inequality~$ao>oc$.
\\
 

\section{Conclusion}
 Until Benedetti's discovery, ovals meant Serlio's ovals, 
 and they were constituted from the `inside'.
 By this, we mean that Serlio took `divine' geometric figures~-- 
 such as circles, squares, equilateral triangles~-- 
 and formed a rhombus, whose edges became the centres of the 4-centre oval.
The  third free parameter  was the radius of one of the circular arcs. All four  ovals by Serlio were constructed in this way.
  
In comparison, Benedetti's ovals were an innovation;
they were  constituted from the `outside'.
Given were the minor and major axes and
an additional parameter, the radius of the larger or smaller
circular arcs. These parameters determine the rhombus with its
centres of the 4-centre oval. The lengths of the axes could be chosen
in a certain ratio, which might guarantee some divineness  of
the construction (for example~$1:2$ or~$3:4$), but this is not necessary~-- 
here, utility could dominate the construction. 
Benedetti's ovals allowed much more flexibility in architectural design.
Of course, this is a theoretical description of the development of ovals
 from  today's art theoretic  point of view.

The historical development, however, took a completely different path. Benedetti was not aware of any architectural needs or applications; we even guess that he did not know about Serlio's oval. It seems that he was a ``true geometer'', 
 driven solely by a problem in Euclidean geometry, 
 namely the construction of the encompassing circle of two circles (Subsection~\ref{Ssec:Bcircles}). Thus, he became intrigued by the problem of constructing an oval for given semiaxes. There is no indication that Benedetti was aware of its mathematical and architectural significance.

His letter to Pizzamano about the oval constructions was published but buried in his~1585 book~\cite{Benedetti-1585} under a heap of letters and other miscellanea.
Nonetheless, about fifteen years later, Clavius unearthed Benedetti's letters and prominently presented them to a wider readership in his main treatises on geometry \cite{Clavius-1604-engl,Euclid-Clavius-1591-engl}, crediting Benedetti with their invention.
In~1623 Beger translated Benedetti's letters into German.
Around the same time, Ardüser, citing only Clavius, provided a precise description of oval~B1 (but not B2) in his German work \textit{Geometriae, Theoricae et Practicae} (1627). Unfortunately, the treatment of oval~B1 concludes with Harsdörffer in 1651, who refers to Benedetti directly but presents only an incomplete description.
Indeed, oval~B1 has not been mentioned in recent overview works on ovals and their history, e.g. \cite{Lopez-Mozo-2011,Mazzotti-2019}.

Oval~B2, on the other hand, was never forgotten but attributed wrongly.
Bosse repeated Benedetti's oval~B2 in~1665, as did Tosca in~1712, but neither of them gave a reference.
This construction was believed to have been invented by Bosse \cite{Mazzotti-2019} or Tosca \cite{Lopez-Mozo-2011}, respectively.
It is plausible that Bordiga's influential biography of Benedetti from~1926, with his unjust, negative verdict on those letters, lastingly kept his readers from studying them in detail.

Summarizing, it appears that Benedetti's ovals~B1 and~B2 constitute the earliest known general constructions for prespecified minor and major axis.
He developed them before~1571, i.e.,  less than~26 years after Serlio and more than 90~years before Bosse.
With their high level of sophistication and mathematical precision, they clearly represent a substantial contribution to the field.


\captionsetup[subfloat]{labelformat=empty}

\ \\[2ex]

\textbf{Acknowledgements}\\[1ex]
We are particularly indebted to Helena Katalin Sz{\'e}pe (University of South Florida) for clarifying biographical details of P.~Pizzamano. Moreover, we are thankful to Clara Silvia Roero (Università di Torino) for hints to the literature on G.\,B.~Benedetti, as well as to Hole Rößler (Herzog August Bibliothek Wolfenbüttel) for pointing us to the work of C.~Clavius.
Also, we thank Angelo~A.\,Mazzotti (Berlin) for fruitful discussions.

\ \\[13ex]

\begin{figure*}[h!]
\begin{center}
 \centering
   \includegraphics[angle=0,width=3.6cm]{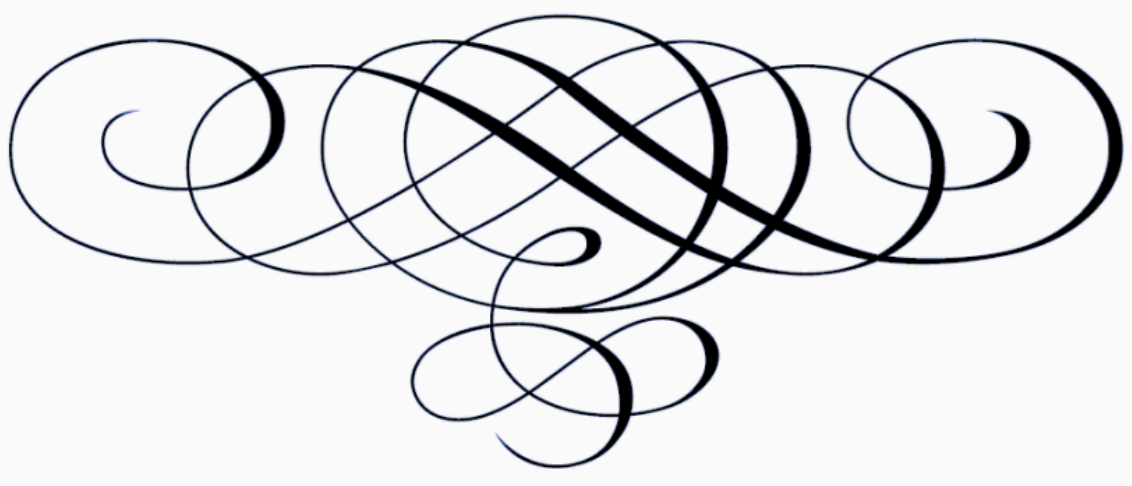}
\end{center}
\end{figure*}

\newpage
\appendix
\section{Appendix: Transcription and translation of Benedetti's letters to 
Pizzamano}

In the following, a slightly modernised and typographically corrected edition of
Benedetti's letters
is presented alongside a translation into English,
the latter favouring intelligibility over literal adherence to the Latin source.
Square brackets contain additions enhancing readability.
\\


\subsection[First letter on constructing a circle encompassing two others]{First letter on constructing a circle encompassing two others\quad \normalfont{(Fig.\,\ref{Fig:Benedetti-262-263-circles})}}\label{App~I}


\  \\
\begin{paracol}{2}
{\noindent\large\bfseries\itshape
Qualiter circulus designari possit alios duos circulos propositos includens
\\[-3pt]}\switchcolumn{\noindent\large\bfseries
How one can describe a circle encompassing two other given circles
\\[-3pt]}\switchcolumn*

{\raggedleft\itshape
Clariss. Petro Pizzamano
\\[12pt]}\switchcolumn{\raggedleft
His Excellency Pietro Pizzamano
\\[12pt]}\switchcolumn*

{\itshape
Superioribus diebus per tuas literas a me quaesivisti, ut modum tibi scribere vellem, quo circulus designari possit circumscribens alios duos propositos circulos.
Qua in re ut tibi satisfaciam quod maxime cupio ita rem accipe.
}\switchcolumn
These past days you asked me through your letter to describe a method to you how one can draw a circle circumscribing two other given circles.
To satisfy you in this matter, what I wish for most, accept the matter this way.
\switchcolumn*

{\itshape
Propositi circuli sint, aut inter se contigui, aut intersecantes vel separati.
Esto primum contiguos esse, qui sint .d.b. et .f.q. quorum .d.b. maior sit et .f.q. minor, eorum vero centra sint .a. et .o. punctum autem contingentiae sit .i.
Nunc subtrahatur .b.a.o.q. per centra eorum ab una circumferentia ad aliam, quae quidem linea transibit per punctum .i. ex .11. tertii Euclidis
deinde a diametro maiori abscindatur .i.e. ad aequalitatem minoris semidiametri,
quo facto sumatur distantia inter .e. et .b. circino mediante factoque centro .o. scindatur, alio circini pede, circumferentia maioris circuli in puncto .u.
a quo si mente concipiemus duas lineas .u.a.d. et .u.o.f. transeuntes per eorum centra .a. et .o. usque ad cicumferentias in punctis .d. et .f. ipsae erunt invicem aequales,
eo quod .e.i. sumpta fuit aequalis .o.f. et .o.u. aequalis .e.b. quare .u.f. aequalis erit .b.i. sed .u.d. etiam aequalis .b.i. ergo .u.d. aequalis erit .u.f. et circulus, cuius .u.d. vel .u.f. erit semidiameter, contiguus erit ipsis propositis circulis ex converso .11. iam dictae.
Idem dico pro circulis se invicem secantibus.
}\switchcolumn
The given circles may either touch each\linebreak
 other, or intersect, or be separated.
First let them,~$db$ and~$fq$, be touching, with~$db$ the larger and~$fq$ the smaller one, their centres~$a$ and~$o$, while~$i$ be the contact point.
Now let the line~$baoq$ be drawn from one circumference through their centres to the\linebreak
 other; it surely will go through the point~$i$ by [Proposition]~11 of Euclid's third [book].
Thereafter cut~$ie$, taken equal to the smaller radius, off the larger diameter.
When that is done, take the distance between~$e$ and~$b$ with the compass and, with~$o$ as centre, cut the point~$u$ off the larger circumference with the other compass leg.
If from there we imagine two lines~$uad$ and~$uof$ running through their centres~$a$ and~$o$ to the circumferences at 
 points~$d$ and~$f$, then they will equal each other,
for~$ei$ was taken equal to~$of$ and~$ou$ equal to~$eb$, whence~$uf$ equals $bi$ but~$ud$ also equals~$bi$, such that~$ud$ equals
$uf$.
The very same I say for intersecting circles.
\switchcolumn*

\newpage
{\itshape
Sed si circuli propositi seiuncti fuerint, sumatur .b.i. diameter maioris, qui fiat semidiameter unius circuli circa centrum .o. et hic circulus vocetur .h.x.
coniungatur deinde semidiameter .o.i. minoris circuli cum semidiametro .a.i. circuli maioris, et ex huiusmodi composita linea, fiat unus semidiameter .a.x. circuli .x.n. concentrici cum maiori,
et a puncto .x. intersectionis horum circulorum (posito quod se invicem intersecent) ducantur per eorum centra .x.a. et .x.o. usque ad ipsorum circumferentias in punctis .d. et .f. duae lineae,
unde habebimus .x.d. aequalem .x.f. eo quod tam in .x.d. quam in .x.f. reperiuntur diametri, et semidiametri amborum circulorum,
facto denique centro .x. unius circuli, cuius semidiameter aequalis sit uni earum .x.d. vel .x.f. solutum erit problema, dicta ratione.
}\switchcolumn
\newpage
But if the given circles were separate, take the diameter~$bi$ of the larger one as the radius of one circle around~$o$; this circle will be called~$hx$.
Thereafter join the smaller circle's radius~$oi$ with the larger circle's radius~$ai$, and let the line thus composed be a radius~$ax$ of the circle~$xn$ with the same centre as the larger one.
And from those circles' intersection point~$x$ (which lies\linebreak
 where they intersect each other) draw two lines~$xa$ and~$xo$ through their centres to their circumferences at points~$d$ and~$f$,\linebreak
whence we have~$xd$ equal to~$xf$ since diameter and radius of the two circles maybe found again both in~$xd$ and in~$xf$.
Now taking~$x$ as the centre of a circle whose radius equals either~$xd$ or~$xf$ will for the\linebreak
 stated reason solve the problem.
\switchcolumn*

{\itshape
Si vero distantia duorum propositorum circulorum tanta fuerit, quod secundi circuli\linebreak
 nequeant se invicem tangere, vel secare, tunc alia via incedendum erit, quae talis est et generalis.
Dividatur tota .q.b. per aequalia in puncto .z. circa quod signentur duo puncta ab ipso aequidistantia .k. et .p.
distantia vero .a.k. facta sit semidiameter esse unius circuli .k.x. circa
 centrum .a. distantia autem .o.p. semidiameter alterius circuli .p.x. circa\linebreak centrum .o.
qui quidem circuli se invicem secent in puncto .x.
a quo cum ductae fuerint .x.a.d. et .x.o.f. per centra dictorum circulorum,
ipse erunt invicem aequales, eo quod cum .b.k. aequalis sit .q.p. igitur .x.d. et .q.p. erunt invicem aequales, sed .f.x. aequalis est .q.p. quare .x.f. aequalis erit x.d.
tunc si .x. centrum fuerit unius circuli, cuius semidiameter sit una
 dictarum, problema solutum erit.
}\switchcolumn
If, however, the distance between the two given circles were so large that the resulting circles could neither touch each other nor intersect, then a different way had to be taken which is the following, and a general one, too.
Let all~$qb$ be divided at point~$z$ in equal parts, from which two points~$k$ and~$p$ are marked at equal distances.
Let in fact the distance~$ak$ be the radius of the circle~$kx$ around the centre~$a$, while the distance~$op$ be the radius of the other circle $px$ around the centre~$o$,
such that these circles intersect at the point~$x$.
From there draw $xad$ and~$xof$ through the centres of said circles;
those themselves equal another, for as~$bk$ equals~$qp$, so~$xd$ and~$qp$ will equal another, but~$fx$ equals~$qp$, whence~$xf$ equals $xd$.
If then~$x$ becomes the centre of one circle whose radius is one of the mentioned ones, the problem will be solved.
\switchcolumn*

{\itshape
Talis etiam solutio commoda erit ad invenien\-dum dictum circulum cuiusvis magnitudinis, dato tamen quod eius diameter, maior sit .b.z. cum in nostra potestate sit accipere puncta .k. et .p. proxima vel remota ab ipso .z. ad libitum.
Unde absque ulla divisione ipsius .q.b. per medium, satis erit signare puncta .k. et .p. duabus distantiis mediantibus .b.k. et .q.p. invicem aequalibus, et etiam propositis.
}\switchcolumn
Such a solution would also be useful in order to obtain said circle at any size, given though that its diameter be larger than~$bz$, which will be in our power by taking 
 points $k$ and~$p$ at will close to or far from~$z$ itself.
Thereby it suffices, without dividing~$qb$ itself in halves, to mark points~$k$ and~$p$ such that their two instrumental distances~$bk$ and $qp$ become equal, still as proposed.
\switchcolumn*
\end{paracol}

\begin{figure*}[t!]
\begin{center}
\subfloat[(a)  \  Oval~B1 for given both semiaxes and radius of large circular arc]
 {
 \centering
   \includegraphics[angle=0,width=6.6cm]{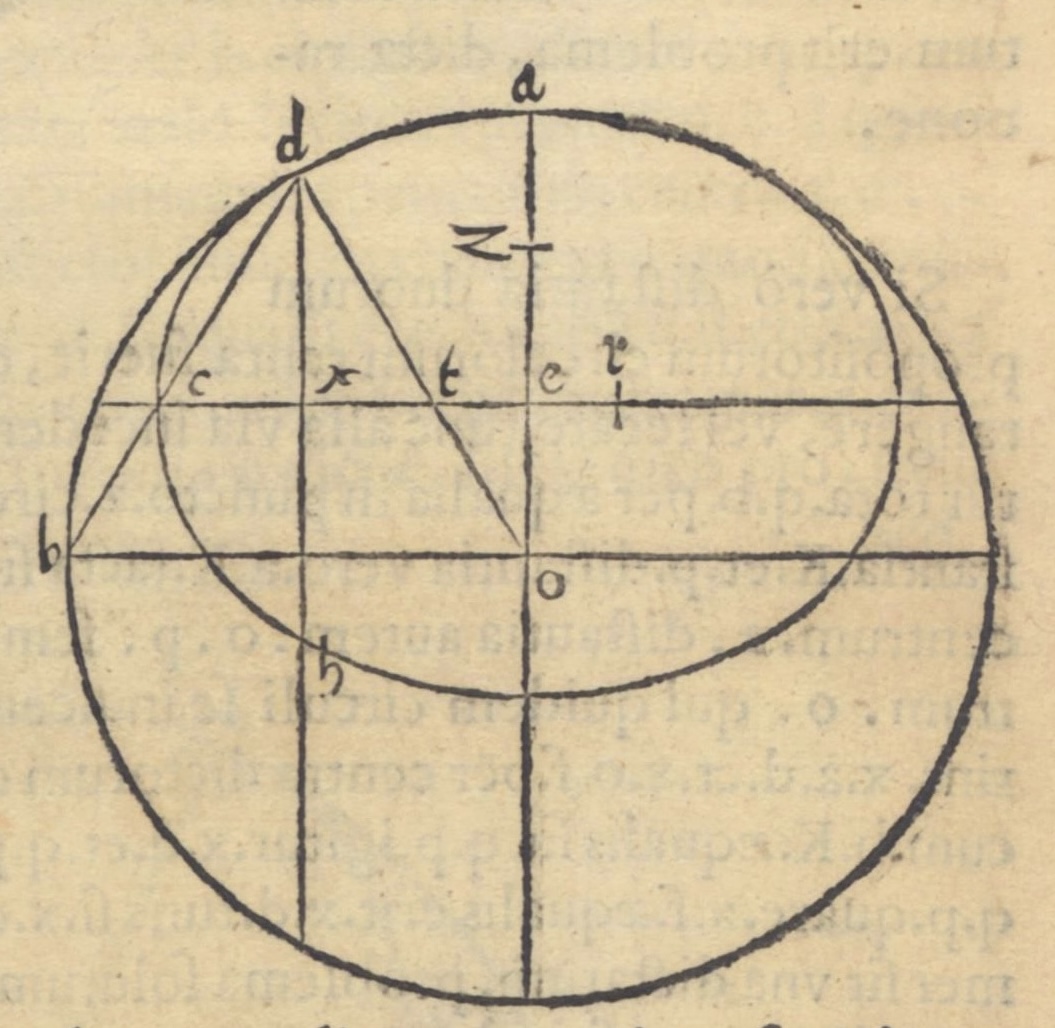}
\label{Abb:Benedetti-1585-264-oben}
 \vspace{5mm}
 }
 \hspace*{6mm}
\subfloat[(b) \  Oval B2 for given both semiaxes and radius of small circular arc]
 {
  \centering
   \includegraphics[angle=0,width=6.6cm]{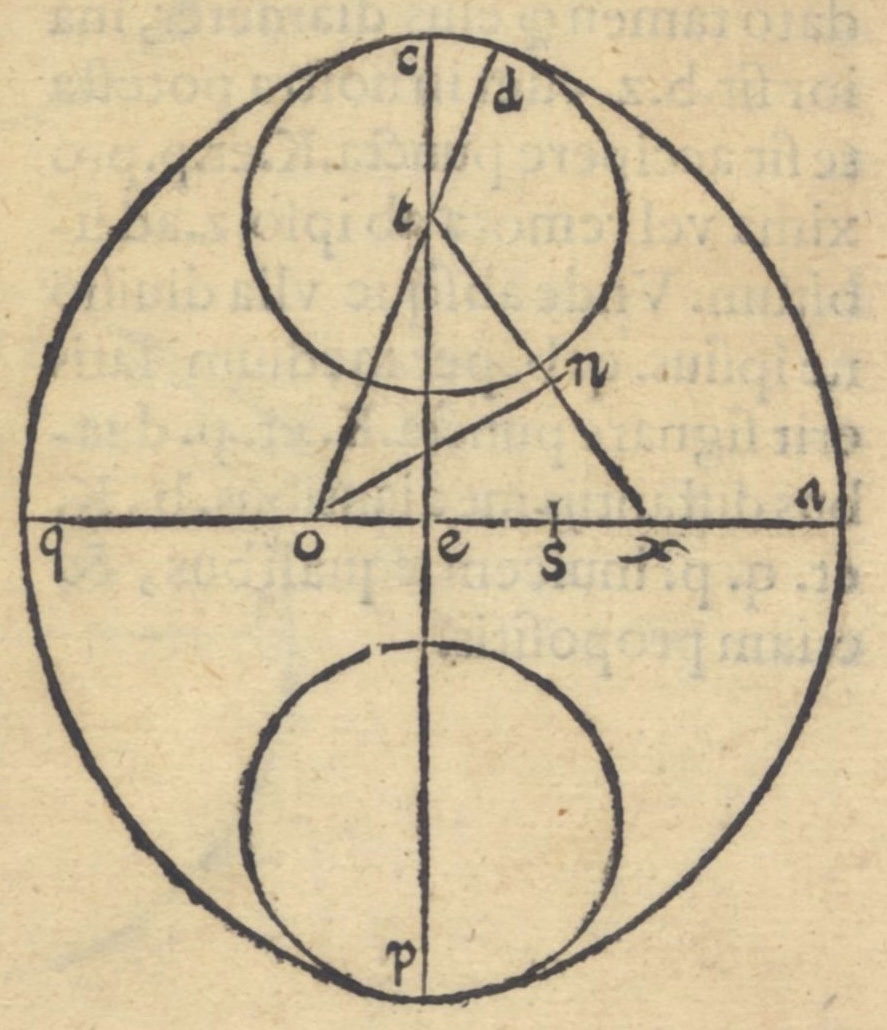}
\label{Abb:Benedetti-1585-264-unten}
}
 \caption{\small  Benedetti's ovals~B1 and~B2 for given semiaxes and one of the circular arcs' radii; details from~\cite[pp.\,264]{Benedetti-1585}. 
} 
 \label{Fig:Benedetti-264-B1-B2}
\end{center}
\end{figure*}

\newpage
\subsection[Second letter on oval constructions]{Second letter on oval constructions\quad \normalfont{(Fig.\,\ref{Fig:Benedetti-264-B1-B2})}}\label{App~II}

\begin{paracol}{2}
{\noindent\large\bfseries\itshape
Figuram superficialem ellipsi similem, ex datis axibus circino mediante delineari posse
\\[-3pt]}\switchcolumn{\noindent\large\bfseries
That it is possible to construct a figure which superficially resembles an ellipse given its axes by means of a compass
\\[-3pt]}\switchcolumn*

{\raggedleft\itshape
Ad eundem
\\[12pt]}\switchcolumn{\raggedleft
To the very same [Clariss. Petro Pizzamano]
\\[12pt]}\switchcolumn*

{\itshape
Figuram superficialem ellipsi similem, ex datis axibus, circino mediante delineare cum volueris, ita facito.
}\switchcolumn
If you wanted to construct a figure which superficially resembles an ellipse given its axes by means of a compass, this is how you ought to do it.
\switchcolumn*

{\itshape
Sit .e.c. semiaxis maior .a.e. vero minor, ad angulum rectum invicem coniuncti, tunc .a.e. producatur usque ad .o.
Itaque .a.o. maior sit quam distantia inter .o. et .c. quae quidem .a.o. posset etiam dari,
describatur postea circulus .a.d.b. circa centrum .o. a quo puncto protrahatur semidiameter .o.b. quae cum .a.o. angulum rectum constituat quae .o.b. erit\linebreak
 aequidistans .e.c. ex .28. primi,
ducatur postea .b.c.d. et .o.t.d. unde angulus .t.c.d. aequalis erit angulo .o.b.d. ex .29. eiusdem.
Ex quinta autem anguli .b. et .d. sunt invicem aequales, quare etiam et anguli .d. et .c. invicem aequales erunt, et ex .6. eiusdem .t.c. aequalis erit .t.d. 
Ducatur postea .d.x.h. perpendicularis lineae .c.e. ita distans sub ipsa .c.e. ut arcus circularis circa .t. delineatus ex semidiametro .t.d. aptus sit eam secare.
Sumpto postea .r. tam distante ab .e. ut .t. reperitur ab ipso .e. et .z. ab .e. ut .o. ab eodem,
ducendo postea duos alios arcus magnitudinis priorum circa centra .r. et .z. habebimus propositum.
}\switchcolumn
Let~$ec$ be the major semiaxis,~$ae$ in fact the minor one, meeting at a right angle, from where~$ae$ will be produced unto~$o$. 
Thus, let~$ao$ be larger than the distance between~$o$ and~$c$ in which way~$ao$ may surely be given;
then describe the circle~$adb$ around its\linebreak
 centre~$o$, from which point the radius~$ob$ is extended forming a right angle with~$ao$, therefore being parallel to~$ec$ by [Proposition]~28 of the first [book of Euclid].
Then draw~$bcd$ and~$otd$ whence by [Proposition]~29 of the same [book] angle~$tcd$ equals angle~$obd$.
By the fifth [proposition] angles~$b$ and~$d$ equal each other, because of which also\linebreak
 angles~$d$ and~$c$ equal each other, and by [Proposition]~6 of the same [book]~$tc$ equals~$td$.
Then~$dxh$ is drawn perpendicular to line~$ce$, and below it until it intersects the circular arc around~$t$ drawn with radius~$td$.
Afterwards taking~$r$ which is to be found as far from~$e$ as~$t$, and~$z$ as far from~$e$ as~$o$ from it,
we will---with two further arcs around~$r$ and~$z$ of the same size as the prior ones---have obtained the proposed.
\switchcolumn*

{\itshape
Sed cum quis voluerit prius arcus minorum circulorum delineare circa maiorem axem, fiant cuiusvis magnitudinis, ut in secunda figura videre est, 
posito tamen quod eorum diameter, minor sit minore axe ipsius figurae, quorum circulorum unus sit .c.d. circa .t. eius centrum,
deinde in axe minori sumatur .a.x. aequalis .c.t. et protrahatur .t.x. quae per aequalia dividatur in puncto .n. a quo postea ducatur .n.o. ad angulos rectos cum .t.x. usque ad intersectionem cum .a.e. in puncto .o. minori axi producta cum oportuerit, 
quod quidem punctum .o. centrum erit arcus .d.a. maioris, eo quod .o.t. aequalis esset .o.x. ex .4. primi Euclidis unde .o.d. aequalis esset .o.a.
et circuli etiam invicem contingentis in puncto .d. ex .11. tertii tam in prima, quam in secunda figura,
sumpto denique puncto .s. tam remoto ab .e. quam .o. reperitur ab eodem, ipsum, centrum erit alterius arcus oppositi,
possemus etiam absque divisione ipsius, .t.x. constituere angulum .x.t.o. aequalem angulo .t.x.o. unde ex .6. primi haberemus .o.t. aequalem .o.x.
}\switchcolumn
But if someone wanted to first draw the smaller circular arcs around the major axis, whatever size they may be, as shown in the second figure,
then given would be that arc's diameter smaller than the smaller axis of the figure, with one of the circles being~$cd$ around its centre~$t$.
Thereafter one takes~$ax$ on the minor axis equal to~$ct$ und draws~$tx$ which is divided equally in the point~$n$,\linebreak
 wherefrom~$no$ is drawn at right angles to~$tx$ until its intersection with~$ae$ in the point~$o$ to which the minor axis should be produced as necessary.
Therefore surely~$o$ will be the centre of the major arc~$da$, because~$ot$\linebreak
 equals~$ox$ by [Proposition]~4 of Euclid's first [book], whence~$od$ equals~$oa$
and the circles touch each other at point~$d$ by [Proposition]~11 of the third [book], both in the first and in the second figure.
Then we take point~$s$, found as far from~$e$ as~$o$, itself as the centre of the other, opposite arc.
Also, without dividing that [line], we could let~$tx$ form the angle~$xto$ equal to angle~$txo$ such that we had~$ot$ equal to~$ox$ by [Proposition]~6 of the first [book].
\switchcolumn*
\end{paracol}

\captionsetup[subfloat]{labelformat=empty}

 \ \\

\addcontentsline{toc}{section}{\hspace*{3.2ex}\textbf{References}}
 
{\small
\bibliographystyle{plain} 

}

 \ \\[12ex]

\textbf{\large Source of figures}\\[1.5ex]
\begin{tabular}{llp{7.5cm}}
Fig.~\ref{Fig:Serlio-ovals}:  & \cite[pp.\,17--19]{Serlio-1545-engl}
&
 Venetia : Sessa, 1551;
München, Bayerische Staatsbibliothek -- Res/2 A.civ. 201-1/5
\\
Fig.~\ref{Abb:Benedetti-Ovals}:  & \qquad ---- & Thomas Hotz
\\
Fig.~\ref{Fig:Benedetti-262-263-circles}: & \cite[pp.\,262--263]{Benedetti-1585}
&
Taurinum, 1585;
München, Bayerische Staatsbibliothek -- 2 Math.u. 0 m
\\
Fig.~\ref{Abb:Cesariano-LXXX}: & \cite[p.\,LXXX verso]{Vitruv-1521-Cesariano-engl}
&
Como : da Ponte, 1521;
Augsburg, Staats- und Stadtbibliothek -- 2 LR 223
 \\
Fig.~\ref{Fig:Clavius}: & \cite[pp.\,218--219]{Clavius-1611-engl}
&
Bamberg, Staatsbibliothek -- RB.Ma.f.1(2\#1
 \\
Fig.~\ref{Abb:Clavius-1591-150}: &\cite[p.\,150]{Euclid-Clavius-1591-engl}
&
Bamberg, Staatsbibliothek -- RB.L.gr.f.1
\\
Fig.~\ref{Abb:Ardueser-1627-S134}: & \cite[p.\,134]{Ardueser-1627}
&
Zürich : Bodmer, 1627;
München, Bayerische Staatsbibliothek -- 4 Math.p. 21
\\
Fig.~\ref{Abb:Bettini_1645_Apiaria-106}: & \cite[Tl.\,1, Apiar.\,III, p.\,106]{Bettini-1645-55-engl}
&
Bononiae : Ferronius, (1648);
Augsburg, Staats- und Stadtbibliothek -- 4 Math 60 -1
\\
Fig.~\ref{Abb:Harsdoerffer3-1653-S180}: & \cite[p.\,180]{Harsdoerffer-1653}
&
Herzog August Bibliothek Wolfenbüttel <\url{http://diglib.hab.de/drucke/224-3-quod/start.htm}>
\\
Fig.~\ref{Abb:Bosse-Oval_1665-T14-fig3}: & \cite[p.\,64 \&\,Plate\,14]{Bosse-1665}
&
A Paris : Chez l'Auteur, M.DC.LXV.
München, Bayerische Staatsbibliothek -- Math.p. 50
 \\
Fig.~\ref{Abb:Tosca-1794-S129-Fig11}: & \cite[Folio~108, Fig.\,11]{Tosca-1727}
&
\url{https://bvpb.mcu.es/es/catalogo_imagenes/grupo.do?path=5571}
\\
Fig.~\ref{Abb:Mazzotti-2019-2nd-Fig31u35}: &  \cite[Fig.\,3.1 and~3.5]{Mazzotti-2019}
&permission by the author
\\
Fig.~\ref{Fig:Benedetti-264-B1-B2}: & \cite[pp.\,264]{Benedetti-1585} 
&
Taurinum, 1585;
München, Bayerische Staatsbibliothek -- 2 Math.u. 0 m
\end{tabular}

\end{document}